# EFFICIENT IMPORTANCE SAMPLING FOR MONTE CARLO EVALUATION OF EXCEEDANCE PROBABILITIES

BY HOCK PENG CHAN[1] AND TZE LEUNG LAI[2]

*National University of Singapore and Stanford University*

Large deviation theory has provided important clues for the choice of importance sampling measures for Monte Carlo evaluation of exceedance probabilities. However, Glasserman and Wang [*Ann. Appl. Probab.* **7** (1997) 731–746] have given examples in which importance sampling measures that are consistent with large deviations can perform much worse than direct Monte Carlo. We address this problem by using certain mixtures of exponentially twisted measures for importance sampling. Their asymptotic optimality is established by using a new class of likelihood ratio martingales and renewal theory.

**1. Introduction.** Importance sampling is a powerful technique to compute the probabilities of rare events by Monte Carlo simulation. For an event occurring with probability $10^{-4}$, one expects the occurrence of 1 event in every 10,000 simulation runs. Therefore, to generate 100 events would require around one million runs for direct Monte Carlo. To simulate a small probability $P(A)$, importance sampling changes the measure $P$ to $Q$ under which $A$ is no longer a rare event and evaluates $P(A)$ by

$$P(A) = \int_A (dP/dQ)\, dQ = E_Q(L\mathbf{1}_A), \tag{1.1}$$

where $L = dP/dQ$ is the likelihood ratio. Whereas $\mathrm{Var}_P(\mathbf{1}_A) = P(A) - (P(A))^2$,

$$\mathrm{Var}_Q(L\mathbf{1}_A) = E_Q(L^2\mathbf{1}_A) - (P(A))^2. \tag{1.2}$$

Therefore, $\mathrm{Var}_Q(L\mathbf{1}_A)$ can be of the order $O((P(A))^2)$ if $L\mathbf{1}_A \le \varepsilon = O(P(A))$, whereas $\mathrm{Var}_P(\mathbf{1}_A) \sim P(A)$ as $P(A) \to 0$.

Received November 2005; revised August 2006.
[1]Supported by research grants from the National University of Singapore.
[2]Supported by the National Science Foundation and the Institute of Mathematical Research at The University of Hong Kong.

In practice, it may be difficult to find $Q$ such that $(dP/dQ)\mathbf{1}_A$ is bounded by $\varepsilon = O(P(A))$. What is needed is $E_Q(L^2\mathbf{1}_A)$ be of the order $O((P(A))^2)$, which is much weaker than $L\mathbf{1}_A$ be bounded by $O(P(A))$. Note that $E_Q L = 1$ whereas $E_Q(L\mathbf{1}_A) = P(A)$ and $A$ is not a rare event under $Q$. Therefore, even though one does not have to choose $Q$ such that $L\mathbf{1}_A$ is bounded by some small number, one has to be careful to avoid the situation where $L\mathbf{1}_A$ is small with a large $Q$-probability but so large with a small $Q$-probability that $E_Q(L^2\mathbf{1}_A)$ is of a larger order of magnitude than $(P(A))^2$. Glasserman and Wang [10] have given examples to show how easily such situations can arise and "how poorly seemingly optimal estimators can perform" when one does not pay attention to avoid such situations. Their paper also gives a brief review of previous work on the choice of $Q$ based on large deviation theory to evaluate exceedance probabilities of random walks, and provides examples for two types of exceedance probabilities which we describe in greater generality below.

Let $\xi, \xi_1, \xi_2, \dots$ be i.i.d. $d$-dimensional random vectors with common distribution $F$ such that $\psi(\theta) := \log(Ee^{\theta'\xi}) < \infty$ for $\|\theta\| < \theta_0$. Let $S_n = \xi_1 + \cdots + \xi_n$, $\mu_0 = E\xi$, $\Theta = \{\theta : \psi(\theta) < \infty\}$, and let $\Lambda$ be the closure of $\nabla\psi(\Theta)$ and $\Lambda^o$ be its interior. Here and in the sequel we use $\nabla\psi$ to denote the gradient vector and $\nabla^2\psi$ the Hessian matrix of second partial derivatives of $\psi$. Then $\nabla\psi$ is a diffeomorphism from $\Theta^o$ onto $\Lambda^o$. Letting $\theta_\mu = (\nabla\psi)^{-1}(\mu)$, define

$$\phi(\mu) = \sup_{\theta \in \Theta}\{\theta'\mu - \psi(\theta)\} = \theta_\mu'\mu - \psi(\theta_\mu), \tag{1.3}$$

which is called the *rate function* in the theory of large deviations. We can embed $F$ in the exponential family $\{F_\theta,\ \theta \in \Theta\}$ with $dF_\theta(x) = e^{\theta'x - \psi(\theta)}\,dF(x)$. Letting $g : \Lambda \to \mathbf{R}$, we consider in Section 2 the exceedance probabilities

$$p_c = P\left\{\max_{n_0 \le n \le n_1} ng(S_n/n) \ge c\right\}, \tag{1.4}$$

$$p_n = P\{g(S_n/n) \ge b\} \qquad \text{with } b > g(\mu_0). \tag{1.5}$$

Let $Q_n$ (or $P_n$) denote the restriction of $Q$ (or $P$) to the $\sigma$-field $\mathcal{F}_n$ generated by $\xi_1, \dots, \xi_n$, and let $P_{\mu,n}$ denote the joint distribution of i.i.d. $\xi_1, \dots, \xi_n$ with common distribution $F_{\theta_\mu}$ and having mean $\mu$. For a stopping time $T$, we also denote the restriction of $Q$ (or $P$, $P_\mu$) to the stopped $\sigma$-field $\mathcal{F}_T$ by $Q_T$ (or $P_T$, $P_{\mu,T}$). In the special case $d = 1$ and $g(x) = x^2$ of (1.5) considered by Glasserman and Wang [10],

$$p_n = P\{|S_n|/n \ge \sqrt{b}\} = P\{|S_n| \ge an\},$$

where $a = \sqrt{b} > |\mu_0|$ and $a \in \Lambda^o$. By large deviation theory, $n^{-1}\log P\{S_n \ge an\} \to -\phi(a)$ and $n^{-1}\log P\{S_n \le -an\} \to -\phi(-a)$. Suppose $\phi(a) < \phi(-a)$. Then $p_n \sim P\{S_n \ge an\}$ and

$$n^{-1}\log L_n \stackrel{P}{\to} -\phi(a) = \lim_{n\to\infty} \log P\{|S_n| \ge an\}, \tag{1.6}$$

where $L_n = dP_n/dP_{a,n}$. Therefore choosing $Q_n = P_{a,n}$ as the importance sampling measure in (1.1) for Monte Carlo computation of $P\{S_n \geq an\}$ is "consistent with large deviations," in the terminology of Glasserman and Wang ([10], page 734), whose Theorem 2 also shows, however, that

$$\lim_{n\to\infty} E_{Q_n}(L_n^2 \mathbf{1}_{\{|S_n|\geq an\}}) = \infty \qquad \text{if } \theta_a + \theta_{-a} > 0. \tag{1.7}$$

Since $\mathrm{Var}_P(\mathbf{1}_{\{|S_n|\geq an\}}) \sim P\{|S_n| \geq an\} = e^{-\{\phi(a)+o(1)\}n}$, (1.7) implies that using the importance sampling measure $Q_n = P_{a,n}$ performs much worse than direct Monte Carlo.

Noting that $A$ has two "minimum rate points" $\pm a$, Glasserman and Wang [10] point out that the preceding difficulty with importance sampling disappears if one uses a mixture $Q_n = \rho P_{a,n} + (1-\rho)P_{-a,n}$ over the minimum rate points $(0 < \rho < 1)$, following an earlier suggestion of Sadowsky and Bucklew [17] who have shown that these mixture-type importance sampling measures are "asymptotically efficient" in the sense that

$$E_{Q_n}(L_n^2 \mathbf{1}_{\{|S_n|\geq an\}}) = e^{-2\{\phi(a)+o(1)\}n}. \tag{1.8}$$

In Section 2 we give a considerably more precise definition of *asymptotic optimality*, replacing the right-hand side of (1.8) by $O(\sqrt{n}p_n^2)$ which we show to be the asymptotically minimal order of the left-hand side over reasonable choices of $Q_n$. More importantly, we provide a much more general way for constructing the asymptotically efficient importance sampling distribution than taking a mixture of $P_{\mu,n}$ over the set of minimum rate points $\mu$, which Sadowsky and Bucklew [17] assume to be a finite set, for general functions $g$ in (1.5).

Glasserman and Wang [10] also consider (1.4) for the special case $d = 2$ and $g(\mu) = \max(\mu_1, \mu_2)$, using $x_j$ to denote the $j$th component of a vector $x$. They assume that $E\xi_{1j} < 0$ for $j = 1, 2$. Setting $n_0 = 1$ and letting $n_1 \to \infty$, this special case of (1.4) reduces to

$$p_c = P\{\max(S_{n,1}, S_{n,2}) \geq c \text{ for some } n \geq 1\} = P\{T_c < \infty\} \sim P\{\tau_c^{(1)} < \infty\}$$

if $\gamma_1 < \gamma_2$, where $\gamma_1$ and $\gamma_2$ are the positive solutions of $\psi(\gamma_1, 0) = 0 = \psi(0, \gamma_2)$ and $\tau_c^{(j)} = \inf\{n : S_{n,j} \geq c\}$, $T_c = \min(\tau_c^{(1)}, \tau_c^{(2)})$. In fact, by Cramér's theorem (cf. [9], page 378), $P\{\tau_c^{(j)} < \infty\} \sim A_j e^{-\gamma_j c}$ (in which $A_j$ is a positive constant not depending on $c$). Glasserman and Wang ([10], Proposition 2) have shown that choosing $Q$ to be the measure under which $\xi_1, \xi_2, \ldots$ are i.i.d. with common distribution $F_{(\gamma_1,0)}$ for Monte Carlo computation of $P\{T_c < \infty\}$ is "consistent with large deviations," in the sense that

$$e^{\gamma_1 c} L_{T_c} \text{ has a nondegenerate limiting distribution as } c \to \infty. \tag{1.9}$$

However, they have also shown that if $\min\{\theta_1 : \psi(\theta_1, \theta_2) = 0 \text{ for some } \theta_2\} > -\gamma_1$, then

$$\lim_{c\to\infty} E_Q(L_{T_c}^2 \mathbf{1}_{\{T_c<\infty\}}) = \infty, \tag{1.10}$$

and therefore this choice of the importance sampling measure $Q$ gives much larger standard error than the direct Monte Carlo estimate of $P\{T_c < \infty\}$, for which $E_P(\mathbf{1}^2_{\{T_c<\infty\}}) = P\{T_c < \infty\} \sim A_1 e^{-\gamma_1 c}$. In Section 2 we resolve this difficulty with importance sampling based on large deviation tilting by using a mixture of the form

$$Q_{T_c\wedge n_1} = \int P_{\mu, T_c\wedge n_1} w_c(\mu)\, d\mu \tag{1.11}$$

for Monte Carlo evaluation of the general boundary crossing probability (1.4). We provide an explicit formula for $w_c(\mu)$ and make use of Theorem 1 of [4] to show that this choice of $Q$ is *asymptotically optimal* in the sense that $E_Q(L_{T_c}^2 \mathbf{1}_{\{T_c<\infty\}})$ attains the asymptotically minimal order of $p_c^2$.

Section 3 generalizes the methods and results of Section 2 to the case where $S_n$ is a Markov random walk, in which $\xi_n$ has distribution $F(\cdot|X_n, X_{n-1})$ depending on a Markov chain $\{X_t\}$. Whereas the methods and results of methods and results of [4] for asymptotic approximations of (1.4) and (1.5) when the increments $\xi_i$ of $S_n$ are i.i.d. provide basic tools for the derivation of the asymptotically optimal importance sampling measure $Q$ in Section 2, the extension to Markov random walks in Section 3 requires new probabilistic ideas. One important idea, given in Section 3.1, is a modification of the usual likelihood ratio martingale to circumvent difficulties with the analysis of eigenfunctions in the Ney–Nummelin [15] formula for likelihood ratios. Section 3.2 develops a new renewal-theoretic approach to the analysis of i.i.d. blocks between regeneration times introduced by Ney and Nummelin [15] for Markov random walks satisfying their minorization condition. Combining these new tools with the results and methods in [5] for the analysis of boundary crossing probabilities, Section 3.3 generalizes (1.11) to Markov random walks. Further refinements of these ideas are used in Section 3.4 for the exceedance probability (1.5).

The complexity due to Markov dependence and nonlinearity in multidimensional settings causes not only analytic difficulties that we resolve in Sections 2 and 3 but also implementation difficulties as the asymptotically optimal importance sampling measure developed in these sections may be difficult to sample directly from. In Section 4 we describe numerical methods to address certain implementation issues and provide numerical examples to illustrate the effectiveness of the methods.

**2. Asymptotically optimal importance sampling measure for Monte Carlo evaluation of exceedance probabilities.** In this section, we derive asymptotically optimal importance sampling measures $Q_c^*$ and $Q_n^*$ for Monte Carlo evaluation of the boundary crossing probability (1.4) and the tail probability (1.5). In particular, it gives an explicit formula (2.1) for a mixing density $w_c(\mu)$ in (1.11) that yields $Q_c^*$. The measure $Q_n^*$ involves a similar mixing density $\widetilde{w}_n(\mu)$ given by (2.13).

2.1. *Boundary crossing probabilities.* Let $T_c = \inf\{n \geq n_0 : ng(S_n/n) \geq c\}$. Then (1.4) can be written as $p_c = P\{T_c \leq n_1\}$. To derive an asymptotically optimal importance sampling measure $Q_c^*$ for Monte Carlo evaluation of $p_c$, we assume the following regularity conditions (A1)–(A5) on $g$ that have been used by Chan and Lai [4] to develop large deviation approximations to $p_c$. Define the rate function $\phi$ by (1.3) and let $\partial\Lambda$ be the boundary of $\Lambda$, $|\cdot|$ denote the determinant of a square matrix, $\Sigma(\mu) = \nabla^2\psi(\theta_\mu)$, and $TM(\mu)$ be the tangent space and $TM^\perp(\mu)$ the normal space of a manifold $M$ at $\mu$.

(A1) There exist $0 < \delta < a < \infty$ and $0 < \varepsilon_0 < a^{-1}$ such that $n_0 \sim \delta c$, $n_1 \sim ac$ and

$$\sup_{a^{-1}-\varepsilon_0 < g(\mu) < \delta^{-1}+\varepsilon_0} g(\mu)/\phi(\mu) = r < \infty.$$

(A2) $M_\varepsilon := \{\mu : a^{-1} - \varepsilon < g(\mu) < \delta^{-1} + \varepsilon \quad \text{and} \quad g(\mu)/\phi(\mu) = r\}$ is a $q$-dimensional oriented manifold for all $0 \leq \varepsilon \leq \varepsilon_0$, where $q \leq d$.

(A3) $\liminf_{\mu\to\partial\Lambda} \phi(\mu) > (\delta r)^{-1}$ and there exists $\varepsilon_1 > 0$ such that $\phi(\mu) > (\delta r)^{-1} + \varepsilon_1$ if $g(\mu) > \delta^{-1} + \varepsilon_0$.

(A4) $g$ is twice continuously differentiable and $\sigma(\{\mu : g(\mu) = \delta^{-1}$ and $g(\mu)/\phi(\mu) = r\}) = 0$, where $\sigma$ is the volume element measure of $M_{\varepsilon_0}$.

(A5) $\inf_{\mu\in M_0} |\nabla_\perp^2\rho(\mu)| > 0$ with $\rho = \phi - g/r$, where $\nabla_\perp^2\rho(\mu) = (\Pi_\mu^\perp)'\nabla^2\rho(\mu)\Pi_\mu^\perp$ and $\Pi_\mu^\perp$ denotes the $d \times (d-q)$ matrix whose column vectors form an orthonormal basis of $TM_0^\perp(\mu)$ in the case $d > q$, and we set $|\nabla_\perp^2\rho(\mu)| = 1$ if $d = q$.

Chan and Lai [4] have given a number of important statistical applications in which (A1)–(A5) are satisfied. In particular, if $g = \phi$, then (A1)–(A5) hold with $r = 1$, $q = d$ and $M_\varepsilon = \{\mu : a^{-1} - \varepsilon < g(\mu) < \delta^{-1} + \varepsilon\}$. The linear function $g(\mu) = r[\theta_{\mu_0}'\mu - \psi(\theta_{\mu_0})]$ also satisfies (A1)–(A5) with $M_\varepsilon = \{\mu_0\}$ and $q = 0$ if $a^{-1} < g(\mu_0) < \delta^{-1}$, but violates (A4) if $g(\mu_0) = \delta^{-1}$. Under (A1)–(A5), let $\Lambda^* = \{\mu \in \Lambda : \phi(\mu) \leq (\delta r)^{-1} + \varepsilon_1$ and $\delta^{-1} + \varepsilon_0 > g(\mu) > a^{-1} - \varepsilon_0\}$ and define

$$\begin{aligned} w_c(\mu) = \beta_c\{&[g(\mu)]^{-d/2}e^{-c\phi(\mu)/g(\mu)}\mathbf{1}_{\{\mu\in\Lambda^*\}} \\ &+ \delta^{d/2}e^{-n_0\phi(\mu)}\mathbf{1}_{\{\phi(\mu)>(\delta r)^{-1}+\varepsilon_1/2\}}\}, \end{aligned} \tag{2.1}$$

where $\beta_c$ is a normalizing constant such that $\int_\Lambda w_c(\mu)\,d\mu = 1$. With this choice of $w_c$, define $Q_c^*$ by the right-hand side of (1.11). The importance sampling method to evaluate $p_c$ by Monte Carlo involves generating $m$ independent samples $(\xi_1^{(i)}, \ldots, \xi_{T_c \wedge n_1}^{(i)})$, $i = 1, \ldots, m$, from $Q_c^*$ so that

$$\widehat{p}_c = m^{-1} \sum_{i=1}^{m} L_c^{(i)} \mathbf{1}_{\{T_c^{(i)} \le n_1\}} \tag{2.2}$$

provides an unbiased estimate of $p_c$, where

$$\begin{aligned} \frac{1}{L_c^{(i)}} &= \frac{dQ_c^*}{dP_{T_c^{(i)} \wedge n_1}} (\xi_1^{(i)}, \ldots, \xi_{T_c^{(i)} \wedge n_1}^{(i)}) \\ &= \int_\Lambda e^{\theta'_\mu S^{(i)}_{T_c \wedge n_1} - (T_c^{(i)} \wedge n_1)\psi(\theta_\mu)} w_c(\mu)\,d\mu. \end{aligned} \tag{2.3}$$

In Section 4, we give details about how to draw the $\xi_t^{(i)}$ from the mixture distribution $Q_c^*$.

To explain the motivation underlying the definition of $w_c(\mu)$, we begin by considering importance sampling to evaluate $P\{S_n/n \in A_n\}$ for a closed bounded set $A_n$ such that $E\xi_1 \notin A_n$. An asymptotically optimal importance density is one that is proportional to $e^{-n\phi(\mu)} \mathbf{1}_{\{\mu \in A_n\}}$. This suggests that to simulate the probability of the event

$$\{T_c \le n_1\} = \{ng(S_n/n) \ge c \text{ for some } n_0 \le n \le n_1\} = \bigcup_{n=n_0}^{n_1} \{ng(S_n/n) \ge c\},$$

it may be optimal to choose an importance density that is proportional to

$$\begin{aligned} &\sup_{n_0 \le n \le n_1} e^{-n\phi(\mu)} \mathbf{1}_{\{ng(\mu) \ge c\}} \\ &\quad = e^{-[c/g(\mu)]\phi(\mu)} \mathbf{1}_{\{c/n_0 \ge g(\mu) \ge c/n_1\}} + e^{-n_0\phi(\mu)} \mathbf{1}_{\{g(\mu) > c/n_0\}}, \end{aligned}$$

in which the supremum on the left-hand side is taken over all real numbers lying between $n_0$ and $n_1$. The formula (2.1) modifies this slightly to facilitate the proof of asymptotic optimality.

We call an importance sampling measure $Q_c$ *asymptotically optimal* for evaluating $p_c$ if

$$E_{Q_c}\left[ \left( \frac{dP_{T_c \wedge n_1}}{dQ_c} \right)^2 \mathbf{1}_{\{T_c \le n_1\}} \right] = O(p_c^2). \tag{2.4}$$

As will be shown in Section 4, there is considerable flexibility in the choice of an asymptotically optimal mixing density. Since $E_{Q_c}[(dP_{T_c \wedge n_1}/dQ_c)\mathbf{1}_{\{T_c \le n_1\}}] = p_c$, the left-hand side of (2.4) is $\ge p_c^2$ by the Cauchy–Schwarz inequality, so the right-hand side of (2.4) indeed gives an asymptotically minimal order to justify the "asymptotic optimality" of (2.4). The following theorem establishes the asymptotic optimality of $Q_c^*$ defined by (1.11) and (2.1).

THEOREM 1. *Assume* (A1)–(A5) *and define* $Q_c^*$ *by* (1.11) *and* (2.1). *Then* $Q_c^*$ *satisfies* (2.4) *and is therefore an asymptotically optimal importance sampling measure.*

PROOF. By considering $g/r$ and $c/r$, we can assume without loss of generality that $r = 1$. We first assume also that $F$ is nonlattice so that Theorem 1 of [4] can be applied, yielding

$$p_c \sim C' c^{q/2} e^{-c} \tag{2.5}$$

for some $C' > 0$. By (2.1) and (A3),

$$\begin{aligned} \beta_c^{-1} &= \int_\Lambda [w_c(\mu)/\beta_c]\, d\mu \\ &\le \delta^{d/2} \int_{\phi(\mu) > \delta^{-1} + \varepsilon_1/2} e^{-n_0 \phi(\mu)}\, d\mu \qquad (2.6) \\ &\quad + (a^{-1} - \varepsilon_0)^{-d/2} e^{-c} \int_{a^{-1} - \varepsilon_0 < g(\mu) < \delta^{-1} + \varepsilon_0} e^{-c\rho(\mu)/g(\mu)}\, d\mu. \end{aligned}$$

Making use of (A5) and arguments similar to those in the proofs of Theorems 1 and 2 of [4], it can be shown that

$$\int_{a^{-1} - \varepsilon_0 < g(\mu) < \delta^{-1} + \varepsilon_0} e^{-c\rho(\mu)/g(\mu)}\, d\mu = O(c^{-(d-q)/2}), \tag{2.7}$$

$$\begin{aligned} \int_{\phi(\mu) > \delta^{-1} + \varepsilon_1/2} e^{-n_0 \phi(\mu)}\, d\mu &= O(n_0^{-(d-1)/2} e^{-n_0(\delta^{-1} + \varepsilon_1/2)}) \\ &= O(c^{-(d-1)/2} e^{-c(1 + \delta\varepsilon_1/3)}). \end{aligned} \tag{2.8}$$

Combining (2.7) and (2.8) with (2.6) yields

$$\beta_c^{-1} = O(c^{(q-d)/2} e^{-c}) \qquad \text{as } c \to \infty. \tag{2.9}$$

Let $B(c; \widehat{\mu}) = \{\mu : \|\mu - \widehat{\mu}\| \le c^{-1/2}\}$. Recalling that $n_0 \sim \delta c$ and $n_1 \sim ac$, we show in the next paragraph that as $c \to \infty$,

$$\beta_c \Big/ \left( \int_{B(c;\hat{\mu})} e^{T[\theta_\mu' \hat{\mu} - \psi(\theta_\mu)]} w_c(\mu)\, d\mu \right) = O(c^{d/2}) \tag{2.10}$$

uniformly in $n_0 \le T \le n_1$ and $Tg(\widehat{\mu}) \ge c$. Let $\bar{\xi}_n = S_n/n$. From (2.9) and (2.10), it follows that

$$\left( \int_\Lambda e^{T_c[\theta_\mu' \bar{\xi}_{T_c} - \psi(\theta_\mu)]} w_c(\mu)\, d\mu \right)^{-2} \mathbf{1}_{\{T_c \le n_1\}} = O(c^q e^{-2c}), \tag{2.11}$$

recalling that $T_c g(\bar{\xi}_{T_c}) \ge c$. In view of (2.3), the desired conclusion (2.4) for $Q_c^*$ follows from (2.5) and (2.11).

To prove (2.10), first consider the case $\inf_{\mu\in B(c;\hat\mu)}\phi(\mu) > \delta^{-1}+\varepsilon_1/2$. Then for $T\ge n_0$,

$$\int_{B(c;\hat\mu)} e^{T[\theta'_\mu\hat\mu-\psi(\theta_\mu)]}[w_c(\mu)/\beta_c]\,d\mu \ge \delta^{d/2}\int_{B(c;\hat\mu)} e^{T[\theta'_\mu\hat\mu-\psi(\theta_\mu)]-T\phi(\mu)}\,d\mu$$
$$= \delta^{d/2}\int_{B(c;\hat\mu)} e^{T\theta'_\mu(\hat\mu-\mu)}\,d\mu,$$

so (2.10) holds. The complementary case $\inf_{\mu\in B(c;\hat\mu)}\phi(\mu)\le\delta^{-1}+\varepsilon_1/2$ implies that there exists $A>0$ such that $\|\widehat\mu\|\le A$, uniformly in $c\ge 1$. Since

$$\sup_{\|\hat\mu\|\le A}\Big[\sup_{\mu\in B(c;\hat\mu)}\phi(\mu)-\inf_{\mu\in B(c;\hat\mu)}\phi(\mu)\Big]\le\varepsilon_1/2$$

for all large $c$, it suffices to consider the case $\sup_{\mu\in B(c;\hat\mu)}\phi(\mu)\le\delta^{-1}+\varepsilon_1$. In this case, for $T\le n_1$ and $Tg(\widehat\mu)\ge c$ with $c$ sufficiently large, $g(\widehat\mu)\ge c/n_1\ge a^{-1}+o(1)$, so $\mu\in\Lambda^*$ for all $\mu\in B(c;\widehat\mu)$. Therefore, letting $\zeta=\inf_{\mu\in\Lambda^*}[g(\mu)]^{-d/2}$,

$$\int_{B(c;\hat\mu)}[w_c(\mu)/\beta_c]\exp\{T[\theta'_\mu\hat\mu-\psi(\theta_\mu)]\}\,d\mu$$
$$\ge\zeta\int_{B(c;\hat\mu)}\exp\{T[\theta'_\mu\hat\mu-\psi(\theta_\mu)]-c\phi(\mu)/g(\mu)\}\,d\mu$$
$$=\zeta\int_{B(c;\hat\mu)}\exp\{T\theta'_\mu(\hat\mu-\mu)$$
$$+[T-c/g(\hat\mu)]\phi(\mu)+c[1/g(\hat\mu)-1/g(\mu)]\phi(\mu)\}\,d\mu$$
$$\ge\zeta e^{-\eta/2}\,\mathrm{Vol}(B(c;\widehat\mu)\cap\{\mu:(\mu-\widehat\mu)'\nabla f(\widehat\mu)\ge 0\}),$$

where $f(\mu)=(T/c)\theta'_\mu(\widehat\mu-\mu)+[1/g(\widehat\mu)-1/g(\mu)]\phi(\mu)$ so that $f(\widehat\mu)=0$, and Taylor's theorem yields $\eta>0$ such that $f(\mu)\ge(\mu-\widehat\mu)'\nabla f(\widehat\mu)-\eta\|\mu-\widehat\mu\|^2/2$ for all $\mu\in B(c;\widehat\mu)$ and large $c$. It then follows that (2.10) also holds when $\sup_{\mu\in B(c;\hat\mu)}\phi(\mu)\le\delta^{-1}+\varepsilon_1$, noting that $\zeta\ge(\delta^{-1}+\varepsilon_1)^{-d/2}$ by (A1), with $r=1$, in this case.

When $F$ is lattice, the preceding arguments can still be used with some minor modifications. In particular, the asymptotic formula (2.5) can be replaced by the weaker result

$$0<\liminf_{c\to\infty}p_c/\{c^{q/2}e^{-c}\}\le\limsup_{c\to\infty}p_c/\{c^{q/2}e^{-c}\}<\infty \tag{2.12}$$

in the lattice case, which suffices to yield (2.4) for $Q_c^*$ from (2.3) and (2.10); see the remark following the proof of Theorem 2. □

2.2. *Tail probabilities of $g(S_n/n)$.* Define

$$\widetilde{w}_n(\mu) = \widetilde{\beta}_n e^{-n\phi(\mu)} \mathbf{1}_{\{g(\mu)\geq b\}}, \qquad \mu \in \Lambda, \tag{2.13}$$

where $\phi$ is the rate function given in (1.3) and $\widetilde{\beta}_n$ is a normalizing constant such that $\int_\Lambda \widetilde{w}_n(\mu)\,d\mu = 1$. Let

$$Q_n^* = \int_\Lambda P_{\mu,n} \widetilde{w}_n(\mu)\,d\mu. \tag{2.14}$$

We propose to use $Q_n^*$ as the importance sampling measure from which $(\xi_1^{(i)}, \ldots, \xi_n^{(i)})$, $i = 1, \ldots, m$, are generated so that

$$\widehat{p}_n = m^{-1} \sum_{i=1}^m L_n^{(i)} \mathbf{1}_{\{g(S_n^{(i)}/n)\geq b\}} \tag{2.15}$$

provides a Monte Carlo estimate of $p_n$, where $S_n^{(i)} = \xi_1^{(i)} + \cdots + \xi_n^{(i)}$ and

$$\frac{1}{L_n^{(i)}} = \frac{dQ_n^*}{dP_n}(\xi_1^{(i)}, \ldots, \xi_n^{(i)}) = \int_\Lambda e^{\theta_\mu' S_n^{(i)} - n\psi(\theta_\mu)} \widetilde{w}_n(\mu)\,d\mu.$$

Note that $\widehat{p}_n$ is an unbiased estimate of $p_n$ with

$$\begin{aligned} \operatorname{Var}(\widehat{p}_n) &= m^{-1} \operatorname{Var}_{Q_n^*}(L_n^{(i)} \mathbf{1}_{\{g(S_n^{(i)}/n)\geq b\}}) \\ &= [E_{Q_n^*}(L_n^2 \mathbf{1}_{\{g(S_n/n)\geq b\}}) - p_n^2]/m. \end{aligned} \tag{2.16}$$

We call an importance sampling measure $Q_n$ *asymptotically optimal* for evaluating the tail probability (1.5) if

$$E_{Q_n}\left[\left(\frac{dP_n}{dQ_n}\right)^2 \mathbf{1}_{\{g(S_n/n)\geq b\}}\right] = O(\sqrt{n} p_n^2). \tag{2.17}$$

Under certain regularity conditions, the following theorem shows that $Q_n^*$ is asymptotically optimal. These regularity conditions are the same as those in Theorem 2 of [4] on large deviation approximations to $P\{g(S_n/n) \geq b\}$, which we restate below using the same notation:

(B1) $g$ is continuous on $\Lambda^o$ and $\inf\{\phi(\mu) : g(\mu) \geq b\} = b/r$ for some $r > 0$.

(B2) $g$ is twice continuously differentiable on $\{\mu \in \Lambda^o : b - \varepsilon_0 < g(\mu) < b + \varepsilon_0\}$ for some $\varepsilon_0 > 0$.

(B3) $\nabla g(\mu) \neq 0$ on $N := \{\mu \in \Lambda^o : g(\mu) = b\}$, and $M := \{\mu \in \Lambda^o : g(\mu) = b, \phi(\mu) = b/r\}$ is a smooth $q$-dimensional manifold (possibly with boundary) for some $0 \leq q \leq d - 1$.

(B4) $\liminf_{\mu\to\partial\Lambda} \phi(\mu) > br^{-1}$ and $\inf_{g(\mu)>b+\delta} \phi(\mu) > br^{-1}$ for every $\delta > 0$.

(B5) $\inf_{\mu\in M}|\Pi'_\mu\{\Sigma^{-1}(\mu)-s\nabla^2 g(\mu)\}\Pi_\mu|>0$ if $d>q+1$, where $s=\|\nabla\phi(\mu)\|/\|\nabla g(\mu)\|$, $e_1(\mu)=\nabla\phi(\mu)/\|\nabla\phi(\mu)\|$, $\{e_1(\mu),e_2(\mu),\dots,e_{d-q}(\mu)\}$ is an orthonormal basis of $TM^\perp(\mu)$ which is a $(d-q)$-dimensional linear space in view of (B3), and $\Pi_\mu$ is the $d\times(d-q-1)$ matrix $(e_2(\mu)\cdots e_{d-q}(\mu))$.

Chan and Lai ([4], pages 1646–1648) have given several important statistical examples in which (B1)–(B5) are satisfied.

Bucklew, Nitinawarat and Wierer [3] have considered an alternative to $\widetilde{w}_n(\mu)\,d\mu$ for the mixing measure in (2.14). Specifically they consider $\widetilde{Q}_n=\int P_{\mu,n}\,d\widetilde{W}(\mu)$, in which unlike (2.1), $\widetilde{W}$ does not depend on $n$ and the distribution of $\xi$ and assigns all its mass to $\{\mu:g(\mu)=b\}$. The price for using these *universal simulation distributions* is that (2.17) has to be replaced by a weaker logarithmic efficiency property

$$E\widehat{p}_n^2=p_n^2e^{o(n)}\qquad\text{as } n\to\infty. \tag{2.18}$$

The following theorem justifies the definition (2.17) of asymptotic optimality by showing that $\sqrt{n}p_n^2$ is the minimal order of magnitude for the left-hand side of (2.17) when $Q_n$ is the joint distribution of i.i.d. $\xi_1,\dots,\xi_n$ with distribution $G$ such that

$$F(A)>0\Rightarrow G(A)>0 \tag{2.19}$$

for any Borel set $A\subset\mathbf{R}^d$, and such that $\lambda(\theta):=\log[\int e^{\theta' x}G(dx)]<\infty$ for all $\|\theta\|\le\theta_1$. More generally, letting $\Gamma=\{\theta:\lambda(\theta)<\infty\}$, $G_\theta$ be the distribution function defined by $dG_\theta(x)=\exp\{\theta' x-\lambda(\theta)\}\,dG(x)$ for $\theta\in\Gamma$, $\widetilde{\theta}_\mu=(\nabla\lambda)^{-1}(\mu)$ and $W_n$ be a distribution function on $\Xi:=\nabla\lambda(\Gamma)$, it considers $Q_n$ of the form

$$Q_n=\int_\Xi Q_{\mu,n}\,dW_n(\mu), \tag{2.20}$$

where $Q_{\mu,n}$ is the joint distribution of i.i.d. $\xi_1,\dots,\xi_n$ with common distribution $G_{\tilde{\theta}_\mu}$.

THEOREM 2. *Assume that $g$ satisfies* (B1)–(B5). *Let $G$ be a distribution function on $\mathbf{R}^d$ satisfying* (2.19) *and such that $\int e^{\theta' x}\,dG(x)<\infty$ for $\theta$ in some neighborhood of the origin. Define $Q_n$ from $G$ via* (2.20), *where $W_n$ is any probability distribution on $\Xi:=\nabla\lambda(\Gamma)$. Then*

$$\liminf_{n\to\infty}E_{Q_n}\left[\left(\frac{dP_n}{dQ_n}\right)^2\mathbf{1}_{\{g(S_n/n)\ge b\}}\right]\Big/(\sqrt{n}p_n^2)>0. \tag{2.21}$$

*Moreover,* (2.17) *holds for $Q_n=Q_n^*$, where $Q_n^*$ is defined by* (2.13) *and* (2.14).

Proof. Dividing $g$ and $b$ by $r$, we assume without loss of generality that $r=1$. To prove that (2.17) holds for $Q=Q_n^*$, let $B(n;\widehat{\mu})=\{\mu:\|\mu-\widehat{\mu}\|\le n^{-1/2}\}$ be a ball of radius $n^{-1/2}$ centered at $\widehat{\mu}$, and we shall show that there exists $\alpha>0$ such that

$$\int_{B(n;\hat{\mu})} e^{n\{\theta'_\mu\hat{\mu}-\psi(\theta_\mu)\}}\widetilde{w}_n(\mu)\,d\mu \ge \alpha n^{-q/2}e^{bn} \qquad \text{whenever } g(\widehat{\mu})\ge b. \tag{2.22}$$

Note that (2.13) yields

$$e^{n\{\theta'_\mu\hat{\mu}-\psi(\theta_\mu)\}}\widetilde{w}_n(\mu)=\widetilde{\beta}_n e^{n\theta'_\mu(\hat{\mu}-\mu)}\mathbf{1}_{\{g(\mu)\ge b\}}. \tag{2.23}$$

We first assume that $F$ is nonlattice so that we can apply Theorem 2 of [4] and its proof to show that for some $C>0$,

$$p_n\sim Cn^{(q-1)/2}e^{-bn}, \tag{2.24}$$

$$\widetilde{\beta}_n^{-1}=\int_{g(\mu)\ge b} e^{-n\phi(\mu)}\,d\mu = O(n^{(q-1-d)/2}e^{-bn}), \tag{2.25}$$

and that there exists $\alpha'>0$ for which

$$\int_{B(n;\hat{\mu})\cap\{\mu:g(\mu)\ge b\}} e^{n\theta'_\mu(\hat{\mu}-\mu)}\,d\mu\ge\alpha' n^{-(d+1)/2} \qquad \text{whenever } g(\widehat{\mu})\ge b. \tag{2.26}$$

Combining (2.23) with (2.25) and (2.26) yields (2.22) for some $\alpha>0$. Let $\bar{\xi}_n=n^{-1}\sum_1^n\xi_i$. Then

$$\begin{aligned} &E_{Q_n^*}\left[\left(\frac{dP_n}{dQ_n^*}\right)^2\mathbf{1}_{\{g(\bar{\xi}_n)\ge b\}}\right]\\ &\quad\le E_{Q_n^*}\left[\left\{\int_{B(n,\bar{\xi}_n)} e^{n(\theta'_\mu\bar{\xi}_n-\psi(\theta_\mu))}\widetilde{w}_n(\mu)\,d\mu\right\}^{-1}\frac{dP_n}{dQ_n^*}\mathbf{1}_{\{g(\bar{\xi}_n)\ge b\}}\right]\\ &\quad\le \alpha^{-1}n^{q/2}e^{-bn}P\{g(S_n/n)\ge b\}\end{aligned} \tag{2.27}$$

by (2.22), noting that $E_{Q_n^*}[(dP_n/dQ_n^*)\mathbf{1}_{\{g(\bar{\xi}_n)\ge b\}}]=P\{g(\bar{\xi}_n)\ge b\}$. From (2.24) and (2.27), it follows that (2.17) holds for $Q_n=Q_n^*$.

To prove that (2.21) holds for $Q_n$ of the form (2.20), we construct neighborhoods $U_n$ of $M$ such that $g(\mu)\ge b$ for $\mu\in U_n$ and

$$\liminf_{n\to\infty} P\{S_n/n\in U_n\}/p_n>0. \tag{2.28}$$

Recall that $e_1(y),\dots,e_{d-q}(y)$ form an orthonormal basis of $TM^\perp(y)$ and that $g=\phi$ on $M$. By (B1) and (B3) with $r=1$, $\phi-g\ge 0$ on $N$ with equality attained on $M$. Hence, for all $y\in M$, $\nabla(\phi-g)(y)\in TN^\perp(y)$. Similarly, $g$ is constant on $N$ and therefore $\nabla g(y)\in TN^\perp(y)$ for all $y\in N$. Since $TN^\perp(y)$ is of dimension 1, it then follows that for every $y\in M$, $\nabla g(y)$ is a scalar multiple of $e_1(y) = \nabla\phi(y)/\|\nabla\phi(y)\|$. For $y \in M$ and $\max_{1\le i\le d-q}|v_i| \le n^{-1/2}$, since $g(y) = b$ and $(\nabla g(y))' \sum_{i=1}^{d-q} v_i e_i(y) = v_1\|\nabla\phi(y)\|/s$, Taylor's expansion yields

$$g\left(y + \sum_{i=1}^{d-q} v_i e_i(y)\right) = b + v_1\|\nabla\phi(y)\|/s + O(v_1^2) - c(v) + o(\|v\|^2), \tag{2.29}$$

where $v = (v_2, \ldots, v_{d-q})'$ and $c(v) = -v'\Pi_y' \nabla^2 g(y)\Pi_y v/2$. Let

$$U_n = \left\{ y + \sum_{i=1}^{d-q} v_i e_i(y) : y \in M, 2n^{-1} \ge v_1 - sc(v)/\|\nabla\phi(y)\| \ge n^{-1}, \max_{2\le i\le d-q}|v_i| \le n^{-1/2} \right\},$$

and note that $g \ge b$ on $U_n$ by (2.29). When $m^{-1}S_m$ has a bounded continuous density $f^{(m)}$ for some $m \ge 1$, the saddlepoint approximation

$$f^{(n)}(\mu) = (1 + o(1))(n/2\pi)^{d/2}|\Sigma(\mu)|^{-1/2}e^{-n\phi(\mu)} \tag{2.30}$$

holds uniformly over compact sets of $\mu$, and we can integrate (2.30) over $U_n$ to obtain

$$P\{S_n/n \in U_n\} = (1 + o(1))(n/2\pi)^{d/2} \int_{U_n} |\Sigma(\mu)|^{-1/2}e^{-n\phi(\mu)}\,d\mu. \tag{2.31}$$

More generally, when $F$ is nonlattice, we can use a tilting argument and a local central limit theorem as in [4], pages 1651–1652, to show that (2.31) still holds. The integral in (2.31) can be evaluated by the same method as that in [4], pages 1650–1653, involving a change of variables for tubular neighborhoods, thereby deriving (2.28) from (2.31) and (2.24).

Let $U_{n,\mu} = \{\sqrt{n}(x - \mu) : x \in U_n\}$ and apply the central limit theorem to conclude that

$$\begin{aligned} &Q_{\mu,n}\{S_n/n \in U_n\} \\ &\quad = Q_{\mu,n}\{n^{-1/2}(S_n - n\mu) \in U_{n,\mu}\} \\ &\quad = \int_{U_{n,\mu}} (2\pi)^{-d/2}|\nabla^2\lambda(\widetilde{\theta}_\mu)|^{-1/2}\exp(-z'\nabla^2\lambda(\widetilde{\theta}_\mu)z/2)\,dz \\ &\qquad + O(n^{-1/2}) \\ &\quad = O(n^{-1/2}) \end{aligned} \tag{2.32}$$

uniformly in $\Xi$. Let $Q_n$ be of the form (2.20). In view of (2.32),

$$Q_n\{S_n/n \in U_n\} = O(n^{-1/2}). \tag{2.33}$$

Letting $\gamma_n = Q_n\{S_n/n \in U_n\}$ and defining the probability measure $\widehat{Q}_n(\cdot) = Q_n(\cdot|S_n/n \in U_n)$, note that $\gamma_n \le \delta n^{-1/2}$ for some $\delta > 0$ by (2.33) and that

$$\begin{aligned} E_{Q_n}\left[\left(\frac{dP_n}{dQ_n}\right)^2 \mathbf{1}_{\{g(S_n/n)\ge b\}}\right] &\ge E_{Q_n}\left[\left(\frac{dP_n}{dQ_n}\right)^2 \mathbf{1}_{\{S_n/n\in U_n\}}\right] \\ &= \gamma_n E_{\hat{Q}_n}[(dP_n/dQ_n)^2] \ge \gamma_n \{E_{\hat{Q}_n}(dP_n/dQ_n)\}^2 \\ &= \gamma_n\{\gamma_n^{-1}P(S_n/n \in U_n)\}^2 \ge \delta^{-1}\sqrt{n}P^2(S_n/n \in U_n). \end{aligned}$$

Therefore (2.21) follows from (2.28).

When $F$ is lattice, we have in place of (2.24),

$$0 < \liminf_{c\to\infty} p_n/\{n^{(q-1)/2}e^{-bn}\} \le \limsup_{c\to\infty} p_n/\{n^{(q-1)/2}e^{-bn}\} < \infty, \tag{2.34}$$

and hence (2.17) follows from (2.25)–(2.27). □

Remark. Suppose $F$ is lattice and let $L_0$ (of full rank $d$) be the minimal lattice of $\xi_1$. In place of (2.30), we now have

$$P\{S_n = u\} = (h_0 + o(1))(2\pi n)^{-d/2}|\Sigma(u/n)|^{-1/2}e^{-n\phi(u/n)}, \tag{2.35}$$

uniformly over compact subsets of $u/n$, with $u \in L_0$, where $h_0 > 0$ is some constant depending only on $L_0$. By summing up (2.35) over $u/n \in U_n$, we obtain

$$P\{S_n/n \in U_n\} = (h_0 + o(1))(2\pi n)^{-d/2} \sum_{u/n\in U_n, u\in L_0} |\Sigma(u/n)|^{-1/2}e^{-n\phi(u/n)},$$

which can be used to replace (2.31) in the preceding argument.

**3. Regeneration, eigenfunctions, eigenmeasures and extension of Theorem 1 to Markov random walks.** Let $\{(X_n, S_n): n = 0, 1, \dots\}$ be a Markov additive process on $\mathcal{X} \times \mathbf{R}^d$ with transition kernel

$$\begin{aligned} P(x, A \times B) &:= P\{(X_1, S_1) \in A \times (B + s)|(X_0, S_0) = (x, s)\} \\ &= P\{(X_1, S_1) \in A \times B|(X_0, S_0) = (x, 0)\}, \end{aligned}$$

for any measurable subset $A \subset \mathcal{X}$, Borel set $B \subset \mathbf{R}^d$ and $s \in \mathbf{R}^d$. We assume that $\{X_n\}$ is aperiodic and irreducible with respect to some maximal irreducibility measure $\varphi$. Let $S_0 = 0$ and define $\xi_n = S_n - S_{n-1}$, so that $S_n = \xi_1 + \cdots + \xi_n$ is a Markov random walk with increments $\xi_i$. We shall assume the minorization condition

$$P(x, A \times B) \ge h(x, B)\nu(A) \tag{3.1}$$

for some probability measure $\nu$ and measure $h(x, \cdot)$ that is positive for all $x$ belonging to a $\varphi$-positive set. Under (3.1) or its variant $P(x, A \times B) \ge$ $h(x)\nu(A \times B)$, Ney and Nummelin [15] have shown that $(X_n, S_n)$ admits a regeneration scheme with i.i.d. inter-regeneration times for an augmented Markov chain, which is called the "split chain." Letting $\tau$ be the first time $(>0)$ to reach the atom of the split chain and assuming that

$$\Omega := \{(\theta, \zeta) : E_\nu e^{\theta' S_\tau - \tau\zeta} < \infty\} \qquad \text{is an open neighborhood of 0,} \tag{3.2}$$

they have shown that for $\theta \in \Theta := \{\theta : (\theta, \zeta) \in \Omega \text{ for some } \zeta\}$, the kernel $\widehat{P}_\theta(x, A) := \int e^{\theta' s} P(x, A \times ds)$ has a maximal simple eigenvalue $e^{\psi(\theta)}$, where $\psi(\theta)$ is the unique solution of the equation $E_\nu e^{\theta' S_\tau - \tau\psi(\theta)} = 1$, with corresponding eigenfunction

$$r(x; \theta) = E_x e^{\theta' S_\tau - \tau\psi(\theta)}. \tag{3.3}$$

Moreover, $\psi(\theta)$ is strictly convex and analytic on $\Theta$ and there exists a full set $F$ [i.e., $\varphi(F^c) = 0$] such that

$$E_x e^{\theta' S_\tau - \tau\zeta} < \infty \qquad \text{for all } x \in F \text{ and } (\theta, \zeta) \in \Omega. \tag{3.4}$$

Therefore, under (3.1) and (3.2), $P$ can be embedded in an exponential family

$$P_\theta(x, dy \times ds) = e^{\theta' s - \psi(\theta)} P(x, dy \times ds) r(y; \theta)/r(x; \theta), \qquad \theta \in \Theta. \tag{3.5}$$

By (3.1) and (3.5), $P_\theta$ satisfies the minorization condition

$$P_\theta(x, A \times B) \geq h_\theta(x, B)\nu_\theta(A) \qquad \text{where } \nu_\theta(dy) = \int_A r(y; \theta)\nu(dy) \tag{3.6}$$

and $h_\theta(x, B) = \int_B h(x, dz) e^{\theta' z - \psi(\theta)}/r(x; \theta)$. Let $\pi(\theta)$ be the stationary distribution under $P_\theta$ and denote $\pi(0)$ simply by $\pi$.

For the special case of i.i.d. $\xi_i$, $e^{\psi(\theta)}$ is the moment generating function $E(e^{\theta' \xi_i})$ and $r(\cdot; \theta) = 1$. Since $r(x; \theta)$ is uniformly positive and bounded under the uniform recurrence condition that there exist $b > a > 0$ and a probability measure $\nu$ on $\mathcal{X} \times \mathbf{R}^d$ for which $a\nu(A \times B) \leq P(x, A \times B) \leq b\nu(A \times B), \forall x \in \mathcal{X}$, and measurable subsets $A$ and $B$ (cf. [11]), it is straightforward to generalize Theorem 1 to uniformly recurrent Markov additive processes. While the uniform recurrence assumption covers the case of finite $\mathcal{X}$, it is too strong for applications to time series and stochastic dynamical systems. Although the same exponential tilting formula (3.5) still holds under the much weaker minorization condition (3.1) than uniform recurrence, $r(X_T; \theta)$ needs no longer be uniformly positive and bounded and its presence in the likelihood ratio statistic $dP_{\theta,T}/dP_{0,T}$ makes the latter intractable. Thus, Ney and Nummelin [15, 16] have to restrict $X_n$ to "s-sets" on which $r(X_n; \theta)$ is within certain bounds when they use (3.5) to analyze large deviation probabilities on $S_n/n$.

To circumvent the intractability of the likelihood ratio statistic, we make use of regeneration times and the representation (3.3) of the eigenfunction to construct a modified likelihood ratio martingale in Section 3.1. We then bound the second moment of the likelihood ratio statistic multiplied by $\mathbf{1}_{\{Tg(S_T/T)\geq c\}}$ by that of the modified likelihood ratio martingale, which we analyze by applying renewal theory to the independent blocks between regeneration times and using an eigenmeasure to bound each of these blocks. Finiteness of the eigenmeasures has been established in Section 3 of [5] under certain "drift conditions" of the type in [14], and we weaken somewhat these conditions in Section 3.1. To highlight the new ideas that are needed for Markov random walks satisfying the minorization condition (3.1), we consider in Section 3.2 the special case $d=1$ and $g(\mu)=\mu$, with $n_0=1$ and $n_1=\infty$, and prove a general theorem (Theorem 4) that yields as corollaries (i) a generalization, to the Markovian setting, of Siegmund's [18] result on asymptotic optimality of $P_{\theta_*,T_c}$ (degenerate mixture over $\theta$) for i.i.d. $\xi_i$, and (ii) a definitive solution of Collamore's [7] closely related problem on simulating ruin probabilities of multidimensional Markov random walks. Theorem 4 is also used to generalize Theorems 1 and 2 to the Markovian setting in Sections 3.3 and 3.4, where comparison with the *dynamic importance sampling* method recently developed by Dupuis and Wang [8] is also given.

3.1. *A modified likelihood ratio martingale.* Let $\mathcal{F}_n$ be the $\sigma$-field generated by $X_0,\dots,X_n,\xi_1,\dots,\xi_n$. Assuming (3.1), Ney and Nummelin ([16], page 596) have shown how a sequence of regeneration times $0<\tau=\tau(1)<\tau(2)<\cdots$ can be constructed with the following three properties: For $k\geq 1$,

$$\tau(k+1)-\tau(k) \text{ are i.i.d. random variables;} \tag{3.7}$$

$$\begin{aligned}&\text{the random blocks } \{X_{\tau(k)},\dots,X_{\tau(k+1)-1},\xi_{\tau(k)+1},\dots,\xi_{\tau(k+1)}\}\\ &\text{are independent;}\end{aligned} \tag{3.8}$$

$$\begin{aligned}&P_x\{X_{\tau(k)}\in A|\mathcal{F}_{\tau(k)-1},\xi_{\tau(k)}\}=\nu(A) \text{ for all } x\in\mathcal{X}\\ &\text{and measurable subsets } A \text{ of } \mathcal{X}.\end{aligned} \tag{3.9}$$

Moreover, for every $n\geq 1$, there exists a measure $h_n(x,\cdot)$ such that

$$P_x\{\tau=n \text{ and } (X_n,\xi_n)\in A\times B\}=\nu(A)h_n(x,B) \qquad \text{for all } x\in\mathcal{X}, \tag{3.10}$$

which is an extension of the regeneration lemma of Athreya and Ney [1] to Markov additive processes.

Set $\tau(0)=0$. Given a stopping time $T$, define the stopping time

$$U=\inf\{u>T: u=\tau(k) \text{ for some } k\geq 1\}. \tag{3.11}$$

For $\theta\in\Theta$, define

$$Z_n(\theta)=\begin{cases}e^{\theta' S_n-n\psi(\theta)}r(X_n;\theta), & \text{if } n<U,\\ e^{\theta' S_U-U\psi(\theta)}, & \text{if } n\geq U.\end{cases} \tag{3.12}$$

Let $\mathcal{G}_n$ be the smallest $\sigma$-field containing $\mathcal{F}_n\cup\sigma\{\tau(k)\mathbf{1}_{\{\tau(k)\leq n\}},k\geq 1\}$.

THEOREM 3. *$Z_n(\theta)$ is a martingale with respect to $\mathcal{G}_n$ under the transition kernel $P$.*

PROOF. For simplicity we shall write $Z_n$ instead of $Z_n(\theta)$. Let $W_n = e^{\theta' S_{n\wedge U} - (n\wedge U)\psi(\theta)}\, r(X_{n\wedge U};\theta)$. Then $W_n$ is a martingale; in fact, (3.5) yields the likelihood ratio martingale

$$\prod_{i=1}^n \{e^{\theta'\xi_i - \psi(\theta)} r(X_i;\theta)/r(X_{i-1};\theta)\} = e^{\theta' S_n - n\psi(\theta)} r(X_n;\theta)/r(x;\theta) \tag{3.13}$$

under $P$. Combining (3.12) with

$$W_n = \begin{cases} e^{\theta' S_n - n\psi(\theta)} r(X_n;\theta), & \text{if } n < U, \\ e^{\theta' S_U - U\psi(\theta)} r(X_U;\theta), & \text{if } n \ge U, \end{cases}$$

and noting that $Z_n = W_n$ on $\{U > n\}$, we obtain

$$\begin{aligned} E[(Z_{n+1} - W_{n+1})\mathbf{1}_{\{U>n+1\}}|\mathcal{G}_n] &= (Z_n - W_n)\mathbf{1}_{\{U>n+1\}} = 0, \\ E[(Z_{n+1} - W_{n+1})\mathbf{1}_{\{U\le n\}}|\mathcal{G}_n] &= (Z_n - W_n)\mathbf{1}_{\{U\le n\}} = Z_n - W_n, \\ E[Z_{n+1}\mathbf{1}_{\{U=n+1\}}|\mathcal{G}_n] &= e^{\theta' S_n - n\psi(\theta)} E_{X_n}[e^{\theta'\xi_1 - \psi(\theta)}\mathbf{1}_{\{\tau=1\}}]\mathbf{1}_{\{T\le n, U>n\}}, \\ E[W_{n+1}\mathbf{1}_{\{U=n+1\}}|\mathcal{G}_n] &= e^{\theta' S_n - n\psi(\theta)} \\ &\quad \times E_{X_n}[e^{\theta'\xi_1 - \psi(\theta)} r(X_1;\theta)\mathbf{1}_{\{\tau=1\}}]\mathbf{1}_{\{T\le n, U>n\}}. \end{aligned}$$

Since $P_x\{\tau = 1 \text{ and } (X_1,\xi_1) \in A \times B\} = \nu(A)h_1(x,B)$ by (3.10) and since $\nu(\mathcal{X}) = 1$ and $\int r(z;\theta)\nu(dz) = E_\nu e^{\theta' S_\tau - \tau\psi(\theta)} = 1$,

$$\begin{aligned} E_x[e^{\theta'\xi_1 - \psi(\theta)}\mathbf{1}_{\{\tau=1\}}] &= \int e^{\theta' z - \psi(\theta)} h_1(x,dz), \\ E_x[e^{\theta'\xi_1 - \psi(\theta)} r(X_1;\theta)\mathbf{1}_{\{\tau=1\}}] &= \left[\int e^{\theta' z - \psi(\theta)} h_1(x,dz)\right]\left[\int r(y;\theta)\nu(dy)\right] \\ &= \int e^{\theta' z - \psi(\theta)} h_1(x,dz). \end{aligned}$$

Therefore, $E[(Z_{n+1} - W_{n+1})\mathbf{1}_{\{U=n+1\}}|\mathcal{G}_n] = 0$. It then follows that $E[(Z_{n+1} - W_{n+1})|\mathcal{G}_n] = Z_n - W_n$. □

The preceding proof shows that for a given stopping time $T$ (in particular the $T_c$ in Section 2.1), we first replace $T$ by the regeneration time $U$ immediately after $T$ and consider the stopped likelihood ratio martingale $W_n$ that replaces $n$ in (3.13) by $n \wedge U$. The modified likelihood ratio martingale (3.12) further replaces $r(X_{n\wedge U};\theta)$ by 1 on the event $\{n \ge U\}$. The reason why this modification helps is that it enables us to bound each of the independent blocks in (3.8) up to the stopping time $U$ by some eigenmeasure of $\mathcal{X}$. For $x \in \mathcal{X}$, define

$$\ell_x(A;\theta,\zeta) = E_x\left[\sum_{n=0}^{\tau-1} e^{\theta' S_n - n\zeta}\mathbf{1}_{\{X_n \in A\}}\right], \tag{3.14}$$

and let $\ell_\nu$ denote $\int \ell_x \, d\nu(x)$. Then $\ell_\nu(\cdot;\theta,\psi(\theta))$ is the left eigenmeasure associated with the eigenvalue $e^{\psi(\theta)}$; see [15, 16]. The finiteness of $\ell_\nu(\mathcal{X};\theta,\psi(\theta))$ and $\ell_x(\mathcal{X};\theta,\psi(\theta))$ has been studied by Chan and Lai ([5], pages 406–409) under certain drift-type conditions. The following lemma considers more generally $\ell_\omega(\theta,\zeta) := \ell_\omega(\mathcal{X};\theta,\zeta)$ instead of requiring $\zeta = \psi(\theta)$, with $\omega = x$ or $\nu$, and can be proved by the same arguments as those used to prove Theorem 4 of [5].

LEMMA 1. *Assume* (3.1) *and* (3.2). *Let* $(\theta,\zeta) \in \Omega$. *Suppose there exist* $0 < \beta < 1$, $a > 0$, *a measurable subset* $C$ *of* $\mathcal{X}$ *with* $\ell_\nu(C;\theta,\zeta) < \infty$ *and* $\ell_x(C;\theta,\zeta) < \infty$ *for all* $x \in \mathcal{X}$, *and a measurable function* $u : \mathcal{X} \to [1,\infty)$ *such that:*

(U1) $E_x[e^{\theta'\xi_1 - \zeta} u(X_1)] \le (1-\beta)u(x)$ *for all* $x \notin C$,
(U2) $\sup_{x\in C} E_x[e^{\theta'\xi_1-\zeta}u(X_1)] \le a$ *and* $\int u(x)\nu(dx) < \infty$.

*Then* $\ell_\nu(\theta,\zeta) < \infty$ *and* $\ell_x(\theta,\zeta) < \infty$ *for all* $x \in \mathcal{X}$.

3.2. *Extension of Siegmund's result on exponential tilting to Markov random walks.* In the case of i.i.d. $\xi_i$ for which $e^{\psi(\theta)}$ is the moment generating function and whose common mean is negative, Siegmund [18] considered the stopping times

$$T_c = \inf\{n \ge 1 : S_n \ge c\}, \qquad T' = \inf\{n \ge 1 : S_n \le -a\}, \tag{3.15}$$

with $0 < a < \infty$, and proposed to use the importance sampling measure $P_{\theta_*}$ for Monte Carlo evaluation of $p_c := P\{T_c < T'\}$, where $\theta_*$ is the unique positive root of $\psi(\theta) = 0$. He also showed that when $P_\theta$ is used as the importance sampling measure, yielding the unbiased estimator

$$\widehat{p}_{\theta,c} := e^{-\theta S_{T_c} + T_c\psi(\theta)}\mathbf{1}_{\{T_c < T'\}}, \tag{3.16}$$

the asymptotically optimal choice of $\theta$ as $c \to \infty$ is $\theta_*$ because $E_{\theta_*}\widehat{p}^2_{\theta_*,c}/E_\theta\widehat{p}^2_{\theta,c} \to 0$ exponentially fast, for all $\theta \ne \theta_*$. Lehtonen and Nyrhinen [12, 13] considered estimation of $p_c$ for $a = \infty$ and showed that the logarithmic efficiency property

$$E\widehat{p}^2_{c,\theta_*} = p_c^2 e^{o(c)} \qquad \text{as } c \to \infty \tag{3.17}$$

holds when the Markov additive process is uniformly recurrent. In this section we make use of the tools developed in Section 3.1 to extend these results to more general Markov random walks and provide a more precise measure of asymptotic efficiency; see Corollary 1. More importantly, we use these tools to prove the following theorem in which the stopping time $T$ need not be of the form (3.15). The theorem, which will be applied in Section 3.3 to generalize Theorem 1 to the Markovian setting, considers the more general $d$-dimensional case and involves the reciprocal $R_n(\theta,\zeta)$ of a modified likelihood ratio statistic which is similar to that in (3.12):

$$R_n(\theta,\zeta)=e^{-\theta' S_n+n\zeta}. \tag{3.18}$$

Let $x_0$ denote the initial state which we assume to belong to the full set $F$ satisfying (3.4).

THEOREM 4. *Assume* (3.1) *and* (3.2). *Let $T$ be a stopping time and define $U$ by* (3.11). *Suppose $(4\theta,4\zeta)$ and $(-2\theta,-2\zeta)$ belong to $\Omega$, $\ell_{x_0}(4\theta,4\zeta)+\ell_\nu(4\theta,4\zeta)<\infty$ and $\theta' E_{\pi(\theta)}\xi_1\neq\zeta$. Then $E_\theta[R_U^2(\theta,\zeta)\mathbf{1}_{\{\theta' S_T-T\zeta\geq c\}}]=O(e^{-2c})$ as $c\to\infty$, where $R_n(\theta,\zeta)$ is defined in* (3.18).

COROLLARY 1. *Let $d=1$ and define $T_c$ by* (3.15) *and*

$$\widehat{p}_{\theta_*,c}=e^{-\theta_* S_{T_c}}[r(x_0;\theta_*)/r(X_{T_c};\theta_*)]\mathbf{1}_{\{T_c<\infty\}}. \tag{3.19}$$

*Assume that $(4\theta_*,0)$ and $(-2\theta_*,0)$ belong to $\Omega$ and that $\ell_{x_0}(4\theta_*,0)+\ell_\nu(4\theta_*,0)<\infty$. Then $E_{\theta_*}\widehat{p}^2_{\theta_*,c}=O(e^{-2\theta_* c})=O(p_c^2)$ and therefore $P_{\theta_*}$ is an asymptotically optimal importance sampling measure.*

PROOF. *Here and throughout the sequel, if the initial state* (*or transition kernel*) *is not specified under the expectation sign, it is assumed to be $x_0$* (*or* $P$). Define $Z_n(\theta)$ by (3.12) with $T=T_c$ in (3.11) and write $Z_n$ instead of $Z_n(\theta_*)$ for simplicity. Since $f(y)=y^{-1}$ is a convex function, $\{Z_n^{-1},\mathcal{G}_n,n\geq 1\}$ is a submartingale under $P$ by Theorem 3. Moreover, since $U>T_c$ by (3.11),

$$Z_{T_c}=e^{\theta_* S_{T_c}}r(X_{T_c};\theta_*). \tag{3.20}$$

Therefore by Jensen's inequality,

$$\begin{aligned}&E\{Z_U^{-1}|T_c<\infty,(X_{T_c},S_{T_c})=(x,s)\}\\&\quad=e^{-\theta_* s}E_x[e^{-\theta_*\tau}]\geq e^{-\theta_* s}(E_x[e^{\theta_*\tau}])^{-1}=e^{-\theta_* s}[r(x;\theta_*)]^{-1}\\&\quad=E\{Z_{T_c}^{-1}|T_c<\infty,(X_{T_c},S_{T_c})=(x,s)\}.\end{aligned}$$

Multiplying the above conditional expectations by $\mathbf{1}_{\{T_c<\infty\}}$ and then taking expectations yields

$$E(Z_U^{-1}\mathbf{1}_{\{T_c<\infty\}})\geq E(Z_{T_c}^{-1}\mathbf{1}_{\{T_c<\infty\}}). \tag{3.21}$$

By (3.9), $X_{\tau(k)}$ is independent of $\{X_1,\dots,X_{\tau(k)-1},\xi_1,\dots,\xi_{\tau(k)}\}$ for all $k\ge 1$, implying that $X_U$ is independent of $(S_U,T_c)$. Therefore,

$$\begin{aligned}E_{\theta_*}(Z_U^{-2}\mathbf{1}_{\{T_c<\infty\}}) &= E_{\theta_*}(e^{-2\theta_* S_U}\mathbf{1}_{\{T_c<\infty\}})\\ &= E[e^{-\theta^* S_U} r(X_U;\theta_*)\mathbf{1}_{\{T_c<\infty\}}]/r(x_0;\theta_*) \qquad (3.22)\\ &= E(Z_U^{-1}\mathbf{1}_{\{T_c<\infty\}})/r(x_0;\theta_*),\end{aligned}$$

noting that $E[r(X_U;\theta)] = \int r(z;\theta)\nu(dz) = 1$. Combining (3.18), (3.21) and (3.22) yields

$$\begin{aligned}E_{\theta_*}[R_U^2(\theta_*,0)\mathbf{1}_{\{\theta_* S_{T_c}\ge\theta_* c\}}] &= E_{\theta_*}(Z_U^{-2}\mathbf{1}_{\{T_c<\infty\}}) = E(Z_U^{-1}\mathbf{1}_{\{T_c<\infty\}})/r(x_0,\theta_*)\\ &\ge E(Z_{T_c}^{-1}\mathbf{1}_{\{T_c<\infty\}})/r(x_0;\theta_*) = E_{\theta_*}(Z_{T_c}^{-2}\mathbf{1}_{\{T_c<\infty\}})\\ &= E_{\theta_*}\widehat{p}^{\,2}_{\theta_*,c}/r^2(x_0;\theta_*),\end{aligned}$$

where the last equality follows from (3.20). In the nonlattice case, $p_c\sim Ae^{-c\theta_*}$ for some constant $A$; see Theorem 3 of [5]. Without the nonlattice assumption, the asymptotic formula can be weakened to $(A_1+o(1))e^{-c\theta_*}\le p_c\le (A_2+o(1))e^{-c\theta_*}$. Since $\theta_* E_{\pi(\theta_*)}\xi_1>0$, Corollary 1 then follows from Theorem 4. $\square$

PROOF OF THEOREM 4. For notational simplicity, denote $R_U(\theta,\zeta)$ by $R_U$. It suffices to show that there exists a constant $B$ such that

$$\begin{aligned}&E_\theta(R_U^2\mathbf{1}_{\{\theta' S_T - T\zeta\ge c\}})\\ &\quad\le e^{-2c}\{[\ell_{x_0}(4\theta,4\zeta)E_{x_0}e^{-2\theta S_\tau+2\tau\zeta}]^{1/2}/r(x_0;\theta) \qquad (3.23)\\ &\qquad + B[\ell_\nu(4\theta,4\zeta)E_\nu e^{-2\theta S_\tau+2\tau\zeta}]^{1/2}\}.\end{aligned}$$

Let $y_k=\theta'[S_{\tau(k)}-S_{\tau(k-1)}]-\zeta[\tau(k)-\tau(k-1)]$ and $\lambda_k=\max_{\tau(k-1)\le n<\tau(k)}\{\theta'[S_n-S_{\tau(k-1)}]-\zeta[n-\tau(k-1)]\}$. By (3.7) and (3.8), the random vectors $(y_k,\lambda_k)$ are i.i.d. for $k\ge 2$. Define the renewal function

$$V(s)=\sum_{k=2}^{\infty}P_\theta\{y_1+\cdots+y_{k-1}\le s\}. \qquad (3.24)$$

Since $E_{\pi(\theta)}\theta'[S_{\tau(k)}-S_{\tau(k-1)}]=\theta'(E_{\pi(\theta)}\xi_1)(E_{\pi(\theta)}\tau)$, $E_\theta y_k\ne 0$ for $k\ge 2$. If $E_\theta y_k>0$, it follows from Blackwell's renewal theorem that there exists a constant $\alpha>0$ such that $V(s+1)-V(s)\le\alpha$ for all $s\in\mathbf{R}$. We can then use this bound in

$$\begin{aligned}&E_\theta(e^{-2\theta' S_U+2U\zeta}\mathbf{1}_{\{\theta' S_T-T\zeta\ge c\}})\\ &\quad=\sum_{k=1}^{\infty}E_\theta(e^{-2(y_1+\cdots+y_k)}\mathbf{1}_{\{\tau(k-1)\le T<\tau(k),\theta' S_T-T\zeta\ge c\}})\end{aligned}$$

$$\begin{aligned}
(3.25)\quad &\le E_\theta(e^{-2y_1}\mathbf{1}_{\{\lambda_1\ge c\}})\\
&\quad + \sum_{k=2}^{\infty}\sum_{s=-\infty}^{\infty} E_\theta(e^{-2(s+y_k)}\mathbf{1}_{\{\lambda_k\ge c-s-1\}})P_\theta\{s<y_1+\cdots+y_{k-1}\le s+1\}\\
&\le e^{-2c}\left[E_\theta(e^{2(\lambda_1-y_1)})+\alpha\sum_{s=-\infty}^{\infty} e^{2(c-s)}E_{\nu_\theta,\theta}(e^{-2y_1}\mathbf{1}_{\{\lambda_1\ge c-s-1\}})\right],
\end{aligned}$$

noting that for $k\ge 2$, $E_\theta(e^{-2y_k}\mathbf{1}_{\{\lambda_k\ge c-s-1\}}) = E_{\nu_\theta,\theta}(e^{-2y_1}\mathbf{1}_{\{\lambda_1\ge c-s-1\}})$ since $X_{\tau(k-1)}$ has distribution $\nu_\theta$; see (3.9) with $\nu$ replaced by $\nu_\theta$. If $E_\theta y_k<0$, then $\sum_{k=2}^{\infty} P_\theta\{s<y_1+\cdots+y_{k-1}\le s+1\}$ is also bounded by $\alpha$ (sufficiently large) for all $s\in\mathbf{R}$, so (3.25) still holds. Note that (3.9) basically says that $X_\tau$ is an "atom" independent of the past history $\{X_0,\ldots,X_{\tau-1},\xi_1,\ldots,\xi_\tau\}$, and therefore in particular is independent of $(y_1,\lambda_1)$. Since $Er(X_\tau;\theta)=\int r(z;\theta)\nu(dz)=1$, it then follows that

$$\begin{aligned}
(3.26)\quad & E_\theta(e^{2(\lambda_1-y_1)})\\
&\quad = E[e^{(2\lambda_1-y_1)}r(X_\tau;\theta)]/r(x_0;\theta)=E(e^{2\lambda_1-y_1})/r(x_0;\theta);
\end{aligned}$$

$$\begin{aligned}
&\sum_{s=-\infty}^{\infty} e^{2(c-s)}E_{\nu_\theta,\theta}(e^{-2y_1}\mathbf{1}_{\{\lambda_1\ge c-s-1\}})\\
&\quad\le \sum_{s=-\infty}^{\infty} e^4\int_{c-s-2}^{c-s-1} e^{2t}E_{\nu_\theta,\theta}(e^{-2y_1}\mathbf{1}_{\{\lambda_1\ge t\}})\,dt\\
(3.27)\quad &\quad = e^4\int_{-\infty}^{\infty} e^{2t}E_{\nu_\theta,\theta}(e^{-2y_1}\mathbf{1}_{\{\lambda_1\ge t\}})\,dt = e^4E_{\nu_\theta,\theta}\left(e^{-2y_1}\int_{-\infty}^{\lambda_1} e^{2t}\,dt\right)\\
&\quad = e^4E_{\nu_\theta,\theta}(e^{2(\lambda_1-y_1)})/2 = e^4E_\nu(e^{2\lambda_1-y_1})/2.
\end{aligned}$$

The last equality of (3.27) follows from

$$E_{\nu_\theta,\theta}(e^{2(\lambda_1-y_1)}) = \int E_x[e^{2\lambda_1-y_1}r(X_\tau;\theta)/r(x;\theta)]r(x;\theta)\,d\nu(x),$$

since $d\nu_\theta(x)=r(x;\theta)\,d\nu(x)$ by (3.6). By the Cauchy–Schwarz inequality and the definition of $\ell_\omega$ in (3.14),

$$(3.28)\quad E_\omega(e^{2\lambda_1-y_1})\le[E_\omega(e^{4\lambda_1})E_\omega(e^{-2y_1})]^{1/2}\le[\ell_\omega(4\theta,4\zeta)E_\omega e^{-2\theta' S_\tau+2\tau\zeta}]^{1/2}$$

for $\omega=x_0$ and $\nu$. From (3.25)–(3.28), (3.23) follows. □

In the case $d>1$, Collamore [7] considered the stopping time

$$(3.29)\quad T_c=\inf\{n: S_n\in cA\}$$

as a generalization of (3.15), where $A \subset \mathbf{R}^d$ and $cA = \{c\mu : \mu \in A\}$. Assume that

(3.30) $A$ is a convex set such that $\partial A$ is a smooth submanifold and $E_\pi \xi_1 \notin cA$ for all $c > 0$.

Then there exist unique $\theta_* \neq 0$ and $\alpha \in \partial A$ such that $\psi(\theta_*) = 0$ and $\theta_*'(\mu - \alpha) \geq 0$ for all $\mu \in A$; see Lemma 3.2 of [7] that proposes to use (3.19), with $\theta_* S_{T_c}$ replaced by $\theta_*' S_{T_c}$, to estimate $p_c = P\{T_c < \infty\}$ in this multidimensional setting. Under certain regularity conditions, Collamore [7] proved the logarithmic efficiency property (3.17). By applying Theorem 4, we can improve (3.17) by providing a much more precise bound on $E_{\theta_*} \widehat{p}_c^2 / p_c^2$, thereby establishing the asymptotic optimality of $P_{\theta_*}$.

COROLLARY 2. *Assume that* (3.30) *holds, that* $(4\theta_*, 0)$ *and* $(-2\theta_*, 0)$ *belong to* $\Omega$ *and that* $\ell_{x_0}(4\theta_*, 0) + \ell_\nu(4\theta_*, 0) < \infty$. *Then*

$$E_{\theta_*} \widehat{p}^2_{\theta_*, c} = O(e^{-2c\theta_*' \alpha}) = O(p_c^2) \qquad \text{as } c \to \infty. \tag{3.31}$$

The derivation of Corollary 2 from Theorem 4 uses the same arguments as those used to prove Corollary 1. In particular, note that

$$p_c = E_{\theta_*}[e^{-\theta_*' S_{T_c}} r(x_0; \theta_*)/r(X_{T_c}; \theta_*)] \sim B e^{-\theta_*' \alpha c}$$

for some constant $B$ in the nonlattice case, as can be shown by a modification of the proof of Theorem 3 of [5]. This asymptotic formula for $p_c$ can be weakened to $(B_1 + o(1)) e^{-\theta_*' \alpha c} \leq p_c \leq (B_2 + o(1)) e^{-\theta_*' \alpha c}$ in the lattice case.

3.3. *Extension of Theorem* 1 *to the Markov setting.* Define $w_c$ by (2.1) and let $Q_c^* = \int P_{\mu, T_c \wedge n_1} w_c(\mu)\, d\mu$, where $P_\mu$ denotes the transition kernel $P_{\theta_\mu}$. The following theorem, whose proof is given in the Appendix, generalizes Theorem 1 to Markov additive processes. It shows that $p_c$ can be estimated efficiently by $\widehat{p}_c = L_c \mathbf{1}_{\{T_c \leq n_1\}}$, where

$$\begin{aligned} \frac{1}{L_c} &= \frac{dQ_c^*}{dP_{T_c \wedge n_1}} (\xi_1, \ldots, \xi_{T_c \wedge n_1}) \\ &= \int_\Lambda \{\exp[\theta_\mu' S_{T_c \wedge n_1} - (T_c \wedge n_1) \psi(\theta_\mu)]\} \\ &\qquad \times r(X_{T_c \wedge n_1}; \theta_\mu) w_c(\mu)\, d\mu / r(x_0; \theta_\mu), \end{aligned} \tag{3.32}$$

and $(\xi_1, \ldots, \xi_{T_c \wedge n_1})$ is generated from $Q_c^*$. Note that the set $\Lambda^*$ in (2.1) has a compact closure under (A1) and (A4).

THEOREM 5. *Assume* (A1)–(A5). *If* $(4\theta_\mu, 4\psi(\theta_\mu))$ *and* $(-2\theta_\mu, -2\psi(\theta_\mu))$ *belong to* $\Omega$ *and* $\ell_{x_0}(4\theta_\mu, 4\psi(\theta_\mu)) + \ell_\nu(4\theta_\mu, 4\psi(\theta_\mu)) < \infty$ *for all* $\mu \in \Lambda^*$, *then* $E_{Q_c^*}(L_c^2 \mathbf{1}_{\{T_c \leq n_1\}}) = O(p_c^2)$, *where* $L_c$ *is defined in* (3.32).

3.4. *Extension of Theorem 2 to Markov additive processes.* First consider the case $d=1$ and $g(\mu)=\mu$. Bucklew, Ney and Sadowsky [2] considered importance sampling for Monte Carlo evaluation of $p_n := P\{S_n/n \ge b\}$ for uniformly recurrent Markov additive processes with $E_\pi \xi_1 < b$. They showed that for all Markov kernels $Q \ne P_b$ satisfying $P \ll Q$, $E_b \widehat{p}_n^2 / E_Q \widehat{p}_{n,Q}^2 \to 0$ exponentially fast as $c \to \infty$, where $\widehat{p}_{n,Q} = (dP/dQ)(\xi_1, \ldots, \xi_n)\mathbf{1}_{\{S_n/n \ge b\}}$ and

$$\widehat{p}_n = e^{-\theta_b S_n + n\psi(\theta_b)}[r(x_0;\theta_b)/r(X_n;\theta_b)]\mathbf{1}_{\{S_n/n\ge b\}}. \tag{3.33}$$

Their proof uses the property that $r(X_n;\theta_b)$ is bounded away from 0 and therefore it suffices to analyze the exponential term $e^{-\theta_b S_n + n\psi(\theta_b)}$. The following theorem, whose proof is given in the Appendix, considers more general Markov additive processes in which the eigenfunctions need not be uniformly positive and show that $P_b$ is still an asymptotically optimal importance measure. It provides a more precise bound on $E_b \widehat{p}_n^2/p_n^2$ than that provided by Bucklew, Ney and Sadowsky [2].

THEOREM 6. *Suppose $d=1$, $g(\mu)=\mu$, and define $\widehat{p}_n$ by (3.33). Assume that $(-2\theta_b, -2\psi(\theta_b))$ and $(4\theta_b, \zeta)$ belong to $\Omega$ for some $\zeta < 4\psi(\theta_b)$ and that $\ell_{x_0}(4\theta_b,\zeta) + \ell_\nu(4\theta_b,\zeta) < \infty$. Then $E_b \widehat{p}_n^2 = O(\sqrt{n} p_n^2)$.*

We next consider the general setting of Theorem 2 and extend it to Markov additive processes. To estimate $p_n := P_{x_0}\{g(S_n/n) \ge b\}$ by Monte Carlo simulations using the importance measure (2.13)–(2.14), the $L_n^{(i)}$ in the estimate (2.15) is given by

$$\begin{aligned}\frac{1}{L_n^{(i)}} &= \frac{dQ_n^*}{dP_n}(\xi_1^{(i)}, \ldots, \xi_n^{(i)}) \\ &= \int_\Lambda e^{\theta_\mu' S_n^{(i)} - n\psi(\theta_\mu)} \widetilde{w}_n(\mu) r(X_n^{(i)};\theta_\mu)\,d\mu / r(x_0;\theta_\mu)\end{aligned} \tag{3.34}$$

in the Markovian setting, where $\theta_\mu$ is the solution of $\nabla\psi(\theta) = \mu$ (see [15], Lemma 3.5, for existence of $\theta_\mu$).

THEOREM 7. *Assume* (B1)–(B5), *and define $\widehat{p}_n$ by* (2.15) *with $L_n^{(i)}$ given by* (3.34). *Assume that for each $\mu$ in a neighborhood of $M$ [see* (B3)], *there exists $\zeta_\mu < 4\psi(\theta_\mu)$ such that $(4\theta_\mu, \zeta_\mu)$ and $(-2\theta_\mu, -2\psi(\theta_\mu))$ belong to $\Omega$ and*

$$\ell_{x_0}(4\theta_\mu, \zeta_\mu) + \ell_\nu(4\theta_\mu, \zeta_\mu) < \infty.$$

*Then $E_{Q_n^*} \widehat{p}_n^2 = O(\sqrt{n} p_n^2)$.*

The proof of Theorem 7 is given in the Appendix. Instead of using the method of mixtures to construct the importance sampling measure, Dupuis and Wang [8] proposed to perform importance sampling via adaptive choice of the tilting parameter at each step to simulate $P\{S_n/n \in A\}$ for uniformly recurrent Markov additive processes. Suppose $(X_k, S_k) = (x, s)$ has been generated. Let $A_k = \{(na - s)/(n - k) : a \in A\}$. Their dynamic importance sampling method chooses $\mu_k \in A_k$ such that $\phi(\mu_k) = \inf\{\phi(a) : a \in A_k\}$ and generates $(X_{k+1}, \xi_{k+1})$ from $P_{\theta_{\mu_k}}(x, \cdot)$. Under certain regularity conditions on $A$, they have established the logarithmic efficiency property (2.18) of the method.

**4. Implementation and examples.** Since $Q_n^*$ is a mixture of $P_{\mu,n}$ with mixing distribution $W_n$ that has density function (2.13) with respect to Lebesgue measure, we can draw the $\xi_j^{(i)}$ from $Q_n^*$ by generating $m$ i.i.d. vectors $(\mu(i), \xi_1^{(i)}, \dots, \xi_n^{(i)})$ as follows: Generate $\mu(i)$ from $W_n$ and then generate $\xi_1^{(i)}, \dots, \xi_n^{(i)}$ from $P_{\mu(i)}$ in the i.i.d. case, and $X_1^{(i)}, \xi_1^{(i)}, \dots, X_n^{(i)}, \xi_n^{(i)}$ from $P_{\mu(i)}$ in the Markov case. These $m$ simulated vectors are used to evaluate $p_n$ by Monte Carlo via (2.15). Likewise, to evaluate $p_c$ by Monte Carlo, we generate $m$ independent trajectories $(\xi_1^{(i)}, \dots, \xi_{T_c^{(i)} \wedge n_1}^{(i)})$ from $P_{\mu(i)}$, where $\mu(i)$ is generated from the distribution with density function $w_c$ given in (2.1). Note that (2.1) and (2.13) involve normalizing constants $\beta_c$ and $\widetilde{\beta}_n$. Instead of using the asymptotically optimal mixture density (2.1), it is often more convenient to use variants thereof that also yield asymptotically optimal importance sampling measures. For example, suppose $v_c(\mu)$ is a density function satisfying

$$\inf_{w_c(\mu)>0} [v_c(\mu)/w_c(\mu)] \geq \varepsilon > 0, \tag{4.1}$$

and let $Q_v$ and $\widehat{p}_v$ denote the corresponding importance sampling measure and associated estimator of $p_c$, respectively. Then $Q_v$ is also an asymptotically optimal importance sampling measure. This property provides us with the flexibility of choosing an importance density $v_c(\mu)$ that does not involve difficult calculation of normalizing constants and such that the likelihood ratio

$$\int_\Lambda e^{\theta'_\mu S_{T\wedge n_1} - (T \wedge n_1)\psi(\theta_\mu)} v_c(\mu) [r(X_T; \theta_\mu)/r(x_0; \theta_\mu)]\, d\mu \tag{4.2}$$

has a closed-form expression or can be easily computed by numerical integration. A statistical application illustrating this point is provided by Chan and Lai ([6], pages 266–268) whose Table 2 shows a large variance reduction over direct Monte Carlo by using a convenient asymptotically optimal importance sampling measure $Q_v$ satisfying (4.1).

Instead of using the mixture of $P_\mu$ with mixing density $w_c$ in (2.1) or $\widetilde{w}_n$ in (2.13), asymptotically optimal importance sampling measures can also be attained by using discrete mixtures of $P_\mu$ whose likelihood ratios do not involve numerical integration. To fix the ideas, first consider the boundary crossing probability $p_c$. Defining

$$K_c(\mu) = \prod_{i=1}^{d} [\mu_i, \mu_i + c^{-1/2}) \qquad \text{for } \mu = (\mu_1, \ldots, \mu_d) \in \mathbf{R}^d, \tag{4.3}$$

and letting $\Lambda_c = \{\mu \in (c^{-1/2}\mathbf{Z})^d : K_c(\mu) \cap \Lambda^* \neq \varnothing\}$, a discrete analogue of (2.1) is the probability mass function

$$\begin{aligned} w_c^*(\mu) = \widehat{\beta}_c\{&[g(\mu)]^{-d/2} e^{-c\phi(\mu)/g(\mu)} \mathbf{1}_{\{\mu \in \Lambda_c\}} \\ &+ \delta^{d/2} e^{-n_0 \phi(\mu)} \mathbf{1}_{\{\phi(\mu) > (\delta r)^{-1} + \varepsilon_1/2\}}\}, \qquad \mu \in (c^{-1/2}\mathbf{Z})^d, \end{aligned} \tag{4.4}$$

where $\widehat{\beta}_c$ is a normalizing constant so that $\sum_{\mu \in (c^{-1/2}\mathbf{Z})^d} w_c^*(\mu) = 1$. The proof of Theorem 5 shows that the theorem still holds if (2.1) is replaced by the probability mass function (4.4). Note that for the special case $d = 1$ and $g(\mu) = \mu$, Corollary 1 only involves a single $P_{\theta_*}$ for the discrete mixture. We next generalize this result to finite mixtures (with support independent of $c$). With $r$ and $\delta$ given in (A1)–(A5), let

$$J(\mu) = \{s \in \Lambda : r[\theta_\mu' s - \psi(\theta_\mu)] \geq \min[\delta^{-1}, g(s)]\}. \tag{4.5}$$

COROLLARY 3. *Assume* (A1)–(A5) *with* $q = 0$, $n_0 = \delta c$ *and* $n_1 = ac$ *and define* $J(\mu)$ *by* (4.5). *Suppose there exists a finite set* $G$ *such that* $\{\mu \in \Lambda : g(\mu) \geq a^{-1}\} \subset \bigcup_{\mu \in G} J(\mu)$. *If* $(4\theta_\mu, 4\psi(\theta_\mu))$ *and* $(-2\theta_\mu, -2\psi(\theta_\mu))$ *belong to* $\Omega$ *and* $\ell_{x_0}(4\theta_\mu, 4\psi(\theta_\mu)) + \ell_\nu(4\theta_\mu, 4\psi(\theta_\mu)) < \infty$ *for every* $\mu \in G$, *then* $\sum_{\mu \in G} \omega_\mu P_\mu$ *is an asymptotically optimal importance sampling measure for any choice of weights* $\omega_\mu$ *such that* $\min_{\mu \in G} \omega_\mu > 0$ *and* $\sum_{\mu \in G} \omega_\mu = 1$.

PROOF. Let $Q = \sum_{\mu \in G} \omega_\mu P_\mu$. Since $n_1 = ac$, $g(S_{T_c}) \geq a^{-1}$ on $\{T_c \leq n_1\}$. Since $\{\mu \in \Lambda : g(\mu) \geq a^{-1}\} \subset \bigcup_{\mu \in G} J(\mu)$ and since $L_{T_c} = dP_{T_c}/dQ_{T_c} \leq \omega_\mu^{-1} dP_{T_c}/dP_{\mu,T_c}$ for every $\mu \in G$, it then follows that

$$\begin{aligned} E_Q[L_{T_c}^2 \mathbf{1}_{\{T_c \leq n_1\}}] &\leq \sum_{\mu \in G} E_Q[L_{T_c}^2 \mathbf{1}_{\{S_{T_c}/T_c \in J(\mu), T_c \leq n_1\}}] \\ &= \sum_{\mu \in G} E[L_{T_c} \mathbf{1}_{\{S_{T_c}/T_c \in J(\mu), T_c \leq n_1\}}] \\ &\leq \sum_{\mu \in G} \omega_\mu^{-1} E_\mu\left[\left(\frac{dP_{T_c}}{dP_{\mu,T_c}}\right)^2 \mathbf{1}_{\{S_{T_c}/T_c \in J(\mu), T_c \leq n_1\}}\right]. \end{aligned} \tag{4.6}$$

Since $T_c \geq n_0 = \delta c$ and $T_c g(S_{T_c}/T_c) \geq c$, it follows from (4.5) that $\{S_{T_c}/T_c \in J(\mu)\} \subset \{\theta_\mu' S_{T_c} - T_c \psi(\theta_\mu) \geq c/r\}$. Therefore, by Theorem 4 [with $\zeta = \psi(\theta_\mu)$] and the proof of Corollary 1,

$$E_\mu[(dP_{T_c}/dP_{\mu,T_c})^2 \mathbf{1}_{\{S_{T_c}/T_c \in J(\mu), T_c \leq n_1\}}] = O(e^{-2c/r}) \tag{4.7}$$

for every $\mu \in G$. Combining (4.6) with (4.7) yields $E_Q[L_{T_c}^2 \mathbf{1}_{\{T_c \le n_1\}}] = O(e^{-2c/r}) = O(p_c^2)$, in view of (2.12) with $q=0$ and with $c$ replaced by the more general form $c/r$. $\square$

Note that Corollary 1 is a special case of Corollary 3 for $d=1$, $g(x)=x$, $\delta=0$, $a=\infty$ and $G=\{\mu_*\}$, where $\mu_* = \psi'(\theta_*)$. In this special case, since $r=\theta_*^{-1}$ and $\psi(\theta_*)=0$, $\{\mu \in \Lambda : g(\mu) \ge a^{-1}\} = \{\mu : \mu \ge 0\} \subset J(\mu_*)$. Finite mixtures of the form $Q = \sum_{\mu \in G} \omega_\mu P_\mu$ are also asymptotically optimal for estimating $p_n = P\{g(S_n/n) \ge b\}$ under conditions similar to those in Corollary 3. Motivated by Glasserman and Wang's [10] example for Monte Carlo evaluation of $P\{|S_n| \ge an\}$, the following corollary of Theorems 4 and 6 considers more general finite mixtures of the form $\sum_{\mu \in G} \omega_{\mu,n} P_{\mu,n}$. Its proof is given in the Appendix.

Corollary 4. *Assume* (B1)–(B5) *for $q=0$. Suppose there exists a finite set $G$ such that $g(\mu) \ge b$ for all $\mu \in G$, $\min_{\mu \in G} \phi(\mu) = b/r$ and $\{\mu \in \Lambda : g(\mu) \ge b\} \subset \bigcup_{\mu \in G} H(\mu)$, where $H(\mu) = \{s \in \Lambda : \theta_\mu'(s-\mu) \ge 0\}$. Assume also that for each $\mu \in G$, $(-2\theta_\mu, -2\psi(\theta_\mu))$ and $(4\theta_\mu, \zeta_\mu)$ belong to $\Omega$ for some $\zeta_\mu < 4\psi(\theta_\mu)$ and $\ell_{x_0}(4\theta_\mu, \zeta_\mu) + \ell_\nu(4\theta_\mu, \zeta_\mu) < \infty$. Then $\sum_{\mu \in G} \omega_{\mu,n} P_{\mu,n}$ is an asymptotically optimal importance sampling measure for any choice of positive weights $\omega_{\mu,n}$ such that $\sum_{\mu \in G} \omega_{\mu,n} = 1$ and*

$$\liminf_{n\to\infty} \omega_{\mu,n}/e^{-2n[\phi(\mu)-b/r]} > 0 \qquad \text{for all } \mu \in G. \tag{4.8}$$

Example 1. For the tail probability $P\{|S_n| \ge an\}$ considered by Glasserman and Wang [10], $\{\mu : |\mu| \ge a\} \subset H(a) \cup H(-a)$, so we can apply Corollary 4 with $G = \{a, -a\}$. Note that their choice of the mixture weights $\omega_{\mu,n} = e^{-n\phi(\mu)}/[e^{-n\phi(a)} + e^{-n\phi(-a)}] = e^{-n\phi(\mu)}/\{(1+o(1))e^{-bn/r}\}$ satisfies (4.8) and that $\min\{\phi(a), \phi(-a)\} = b/r$. We study the performance of this importance sampling measure in a more general example of a Markov additive process in which the underlying Markov chain $\{X_n\}_{n\ge 0}$ has state space $\{1,2,3\}$ and transition matrix $(p_{xy})_{x,y\in\mathcal{X}}$ such that $p_{xx}=0.5$ for every $x$, $p_{12}=p_{23}=p_{31}=0.3$, $p_{13}=p_{21}=p_{32}=0.2$. Letting $\xi_i = X_i$ so that $S_n = X_1 + \cdots + X_n$, consider the Monte Carlo evaluation of $P_\pi\{S_n/n \ge 2.7 \text{ or } S_n/n \le 1.5\}$, where $\pi$ is the stationary distribution with $\pi(1)=\pi(2)=\pi(3)=1/3$; this corresponds to Corollary 4 with $g(\mu) = (\mu - 2.1)^2$ and $\sqrt{b} = 0.6$. Table 1 compares direct Monte Carlo evaluation of this probability with two importance sampling procedures, the first using $Q_n = P_{1.5,n}$ that tilts to the minimum rate point $\mu = 1.5$, and the second using the mixture

$$Q_n = \omega_n P_{1.5,n} + (1-\omega_n) P_{2.7,n} \qquad \text{with } \omega_n = e^{-n\phi(1.5)}/(e^{-n\phi(1.5)} + e^{-n\phi(2.7)}), \tag{4.9}$$

as advocated by Glasserman and Wang **(year?)**. The eigenvalues and eigenvectors used to define $P_{1.5,n}$ and $P_{2.7,n}$ are given by

$$\begin{aligned}
\theta_{1.5} &= -0.507, & e^{\psi(\theta_{1.5})} &= 0.688,\\
\theta_{2.7} &= 0.815, & e^{\psi(\theta_{2.7})} &= 3.11,\\
r(1;\theta_{1.5}) &= 1.20, & r(1;\theta_{2.7}) &= 0.747,\\
r(2;\theta_{1.5}) &= 0.88, & r(2;\theta_{2.7}) &= 1.02,\\
r(3;\theta_{1.5}) &= 0.92, & r(3;\theta_{2.7}) &= 1.23.
\end{aligned}$$

Moreover, $\phi(1.5) = 0.120$ and $\phi(2.7) = 0.251$ are used to evaluate $\omega_n$. Each result in Table 1 is a Monte Carlo estimate based on $m =$10,000 simulation runs, and the standard error of the estimate is given in parentheses. Table 1 shows that direct Monte Carlo has much larger standard error than importance sampling with (4.9) and that it is unable to provide a meaningful estimate when the probabilities are smaller than $10^{-4}$. The minimum rate point method ($Q_n = P_{1.5,n}$) has much larger standard error than (4.9) for $n = 10$ and tends to underestimate the true probability for $n = 20$.

EXAMPLE 2. Let $\varepsilon_1, \varepsilon_2, \ldots$ be i.i.d. three-dimensional standard normal vectors and $1 \le n_0 \le n_1 < \infty$. Consider the *regime-switching* Gaussian random walk $S_n = \sum_{i=1}^n \xi_i$, with $\xi_i = (X_i - 2, 0, 0)' + \varepsilon_i$, where $\{X_n\}_{n\ge 0}$ is the Markov chain in Example 1. Let $T_c = \inf\{n \ge n_0 : \|S_n\|^2 \ge cn\}$, which corresponds to the case $g(\mu) = \|\mu\|^2$ in Theorem 5. To compute $p_c = P_\pi\{T_c \le n_1\}$ via importance sampling, where $\pi = (1/3, 1/3, 1/3)'$ is the stationary distribution of $X_0$, we use a slight modification of (4.4) to define the discrete mixture density function by

$$\begin{aligned}
\bar{w}_c(\mu) = \bar{\beta}_c\{&[g(\mu)]^{-3/2} e^{-c\phi(\mu)/g(\mu)} \mathbf{1}_{\{c/n_1 \le g(\mu) \le c/n_0\}} \\
&+ (c/n_0)^{3/2} e^{-n_0\phi(\mu)} \mathbf{1}_{\{c/n_0 < g(\mu) \le b\}}\},
\end{aligned} \tag{4.10}$$

TABLE 1

| $n$ | Direct Monte Carlo | Importance sampling: $Q_n = P_{1.5,n}$ | Importance sampling: $Q_n$ given by (4.9) |
|---|---|---|---|
| 10 | $1.04(0.03) \times 10^{-1}$ | $1.1(0.3) \times 10^{-1}$ | $1.04(0.01) \times 10^{-1}$ |
| 20 | $2.0(0.1) \times 10^{-2}$ | $1.91(0.03) \times 10^{-2}$ | $2.02(0.03) \times 10^{-2}$ |
| 40 | $1.1(0.3) \times 10^{-3}$ | $1.25(0.02) \times 10^{-3}$ | $1.25(0.02) \times 10^{-3}$ |
| 60 | $2(1) \times 10^{-4}$ | $0.93(0.02) \times 10^{-4}$ | $0.96(0.02) \times 10^{-4}$ |
| 80 | 0 | $7.4(0.2) \times 10^{-6}$ | $7.4(0.2) \times 10^{-6}$ |
| 100 | 0 | $5.9(0.1) \times 10^{-7}$ | $5.9(0.1) \times 10^{-7}$ |

for $\mu \in (c^{-1/2}\mathbf{Z})^3$, with $b \geq c/n_0$ and normalizing constant $\bar{\beta}_c$ such that $\sum_{\mu\in(c^{-1/2}\mathbf{Z})^3} \bar{w}_c(\mu) = 1$. Note that

$$\widehat{P}_\theta := (E_x[e^{\theta'\xi_1}\mathbf{1}_{\{X_1=y\}}])_{1\leq x,y\leq 3} = e^{\|\theta\|^2/2}\widetilde{P}_a,$$

where $a$ is the first component of $\theta$ and

$$\widetilde{P}_a = \begin{pmatrix} 0.5e^{-a} & 0.3 & 0.2e^a \\ 0.2e^{-a} & 0.5 & 0.3e^a \\ 0.3e^{-a} & 0.2 & 0.5e^a \end{pmatrix}.$$

Let $\lambda(a)$ be the largest log-eigenvalue of $\widetilde{P}_a$, with associated eigenvector $r(\cdot;a)$. Then $\psi(\theta_\mu) = \lambda(a_\mu) + \|\theta_\mu\|^2/2$, where we use superscripts to denote the components of the vectors $\mu = (\mu^1, \mu^2, \mu^3)'$ and $\theta_\mu = (a_\mu, \theta_\mu^2, \theta_\mu^3)'$. Let $\dot{\lambda}$ denote the derivative of $\lambda$. Since $\mu = \nabla\psi(\theta_\mu) = (\dot{\lambda}(a_\mu), 0, 0)' + \theta_\mu$, $\theta_\mu^2 = \mu^2$ and $\theta_\mu^3 = \mu^3$. Moreover, $\lambda$ is convex and therefore we can use the bisection method to solve the equation $\dot{\lambda}(a) + a = \mu^1$ for $a = a_\mu$. We first generate $\mu$ from the mixture density function (4.10) and then use the measure $P_{\theta_\mu}$ to generate $\xi_i = (X_i - 2, 0, 0)' + \theta_\mu + \varepsilon_i$ so that $\{X_n\}_{n\geq 0}$ has the transition probabilities

$$P_{a_\mu}(x,y) = \widetilde{P}_{a_\mu}(x,y)e^{-\lambda(a_\mu)}r(y;a_\mu)/r(x;a_\mu).$$

Table 2 gives Monte Carlo estimates of $p_c$ for eight choices of $(c, n_0, n_1)$, based on $m = 10{,}000$ runs for each entry, in which the standard error is shown in parentheses. We compare direct Monte Carlo with importance sampling using the mixture density function (4.10) in which $b = 7$. The results show the effectiveness of using (4.10) to compute probabilities of order as small as $10^{-7}$, and that direct Monte Carlo becomes unreliable even for probabilities of order $10^{-4}$. Although extra time is used by importance sampling to compute the likelihood ratio $L_{T_c}$, direct Monte Carlo and importance sampling have similar simulation times because direct Monte Carlo has to generate $X_n$ until $n = n_1$ for most runs whereas importance sampling can stop at $T_c < n_1$.

TABLE 2

| $c$ | $n_0$ | $n_1$ | **Direct Monte Carlo** | **Importance sampling with weights (4.10)** |
|---|---|---|---|---|
| 20 | 5 | 50 | $3.0(0.2) \times 10^{-2}$ | $3.19(0.05) \times 10^{-2}$ |
| 25 | 5 | 50 | $9(1) \times 10^{-3}$ | $8.57(0.08) \times 10^{-3}$ |
| 30 | 5 | 50 | $3.1(0.6) \times 10^{-3}$ | $2.75(0.03) \times 10^{-3}$ |
| 35 | 5 | 50 | $5(2) \times 10^{-4}$ | $5.58(0.07) \times 10^{-4}$ |
| 40 | 10 | 100 | $4(2) \times 10^{-4}$ | $7.3(0.3) \times 10^{-4}$ |
| 50 | 10 | 100 | 0 | $3.37(0.09) \times 10^{-5}$ |
| 60 | 10 | 100 | 0 | $1.82(0.04) \times 10^{-6}$ |
| 70 | 10 | 100 | 0 | $9.2(0.2) \times 10^{-8}$ |

## APPENDIX

PROOF OF THEOREM 5. Assume without loss of generality $r = 1$. First note that (2.9) is still valid. Moreover, (2.12) still holds, as can be shown by arguments similar to the proof of Theorem 6 of [5]. Define $\Gamma_c = \{\mu : \phi(\mu) \le \delta^{-1} + \varepsilon_1\} \cap (c^{-1/2}\mathbf{Z})^d$ and $A_c(\mu^*) = \{T_c \le n_1, S_{T_c}/T_c \in K_c(\mu^*)\}$, where $K_c(\mu)$ is defined in (4.3). We next apply Theorem 4 to show that uniformly in $\mu^* \in \Gamma_c$ with $K_c(\mu^*) \cap \Lambda^* \neq \varnothing$,

$$\begin{aligned} &E_{Q_c^*}(L_c^2 \mathbf{1}_{A_c(\mu^*)}) \\ &\qquad = O\biggl(\beta_c^{-1} c^{d/2} \exp\biggl\{c\biggl[-2 \inf_{\mu \in K_c(\mu^*)} \phi(\mu)/g(\mu) \\ &\qquad\qquad + \sup_{\mu \in K_c(\mu^*)} \phi(\mu)/g(\mu)\biggr]\biggr\}\biggr). \end{aligned} \tag{A.1}$$

Define $Z_n(\theta)$ as in (3.12). Let $\eta_{c,\mu^*} = \int_{K_c(\mu^*)} w_c(\mu)\, d\mu$ and $\widetilde{w}_c(\mu) = w_c(\mu)/\eta_{c,\mu^*}$. Then

$$\begin{aligned} E_{Q_c^*}(L_c^2 \mathbf{1}_{A_c(\mu^*)}) &= E(L_c \mathbf{1}_{A_c(\mu^*)}) \\ &\le E\biggl[\biggl(\int_{K_c(\mu^*)} Z_{T_c}(\theta_\mu) w_c(\mu)\, d\mu\biggr)^{-1} \mathbf{1}_{A_c(\mu^*)}\biggr] \\ &= \eta_{c,\mu^*}^{-1} E\biggl[\biggl(\int_{K_c(\mu^*)} Z_{T_c}(\theta_\mu) \widetilde{w}_c(\mu)\, d\mu\biggr)^{-1} \mathbf{1}_{A_c(\mu^*)}\biggr]. \end{aligned} \tag{A.2}$$

Since $\int_{K_c(\mu^*)} \widetilde{w}_c(\mu)\, d\mu = 1$, Jensen's inequality yields

$$\biggl(\int_{K_c(\mu^*)} Z_{T_c}(\theta_\mu) \widetilde{w}_c(\mu)\, d\mu\biggr)^{-1} \le \int_{K_c(\mu^*)} Z_{T_c}^{-1}(\theta_\mu) \widetilde{w}_c(\mu)\, d\mu.$$

Putting this in (A.2), we obtain

$$\begin{aligned} E_{Q_c^*}(L_c^2 \mathbf{1}_{A_c(\mu^*)}) &\le \eta_{c,\mu^*}^{-1} \int_{K_c(\mu^*)} E[Z_{T_c}^{-1}(\theta_\mu) \mathbf{1}_{A_c(\mu^*)}] \widetilde{w}_c(\mu)\, d\mu \\ &\le \eta_{c,\mu^*}^{-1} \sup_{\mu \in K_c(\mu^*)} E[Z_{T_c}^{-1}(\theta_\mu) \mathbf{1}_{A_c(\mu^*)}]. \end{aligned} \tag{A.3}$$

Since the function $h_\mu(\theta) := \theta'\mu - \psi(\theta)$ is maximized at $\theta = \theta_\mu$ with maximum value $\phi(\mu)$, $\nabla h_\mu(\theta_\mu) = 0$ and there exists $\alpha_1 > 0$ such that $h_{S_{T_c}/T_c}(\theta_\mu) \ge$ $\phi(S_{T_c}/T_c) - \alpha_1 c^{-1}$ if $\mu$ and $S_{T_c}/T_c$ belong to $K_c(\mu^*)$. Therefore on $A_c(\mu^*)$ and for $\mu \in K_c(\mu^*)$,

$$\begin{aligned} &\theta_\mu' S_{T_c} - T_c\psi(\theta_\mu) \\ &\quad = T_c h_{S_{T_c}/T_c}(\theta_\mu) \geq T_c\phi(S_{T_c}/T_c) - T_c\alpha_1 c^{-1} \\ &\quad \geq c\phi(S_{T_c}/T_c)/g(S_{T_c}/T_c) - \alpha_2 \\ &\quad \geq c \inf_{\mu\in K_c(\mu^*)}[\phi(\mu)/g(\mu)] - \alpha_2. \end{aligned} \tag{A.4}$$

In view of (A.4), application of Theorem 4 with $\zeta = \psi(\theta_\mu)$ and Jensen's inequality as in the proof of Corollary 1 then shows that for $\mu \in K_c(\mu^*)$,

$$\begin{aligned} &E[Z_{T_c}^{-1}(\theta_\mu)\mathbf{1}_{A_c(\mu^*)}] \\ &\quad \leq E[Z_U^{-1}(\theta_\mu)\mathbf{1}_{A_c(\mu^*)}] \\ &\quad = r(x_0;\theta_\mu)E_{\theta_\mu}[Z_U^{-2}(\theta_\mu)\mathbf{1}_{A_c(\mu^*)}] \\ &\quad = O\left(\exp\left\{-2c \inf_{\mu\in K_c(\mu^*)}[\phi(\mu)/g(\mu)]\right\}\right). \end{aligned} \tag{A.5}$$

From (2.1), it follows that uniformly in $\mu^* \in \Gamma_c$ with $K_c(\mu^*) \cap \Lambda^* \neq \varnothing$,

$$\eta_{c,\mu^*}^{-1} = O\left(\beta_c^{-1} c^{d/2} \exp\left\{c \sup_{\mu\in K_c(\mu^*)}[\phi(\mu)/g(\mu)]\right\}\right). \tag{A.6}$$

Combining (A.3) with (A.5) and (A.6) yields (A.1).

Making use of (A1)–(A5), we can use geometric integration as in [4], page 1651, to show that

$$\begin{aligned} \sum_{\mu^*\in\Gamma_c:K_c(\mu^*)\cap\Lambda^*\neq\varnothing} &\exp\left\{c\left[-2\inf_{\mu\in K_c(\mu^*)}\left(\frac{\phi(\mu)}{g(\mu)} - 1\right)\right.\right. \\ &\left.\left. + \sup_{\mu\in K_c(\mu^*)}\left(\frac{\phi(\mu)}{g(\mu)} - 1\right)\right]\right\} = O(c^{q/2}). \end{aligned}$$

Combining this with (A.1), in which $\beta_c^{-1} = O(c^{(q-d)/2}e^{-c})$, then yields

$$\sum_{\mu^*\in\Gamma_c:K_c(\mu^*)\cap\Lambda^*\neq\varnothing} E_{Q_c^*}(L_c^2\mathbf{1}_{A_c(\mu^*)}) = O(c^q e^{-2c}) = O(p_c^2) \tag{A.7}$$

by (2.12). Moreover, the proof of Lemma 2 in [5], pages 418–419, can be used to show that

$$E_{Q_c^*}(L_c^2\mathbf{1}_{\{T_c\leq n_1, S_{T_c}/T_c\notin\Lambda^*\}}) = o(e^{-2c}). \tag{A.8}$$

From (A.7) and (A.8), the desired conclusion follows. □

PROOF OF THEOREM 6. Let $U$ be the stopping time (3.11) associated with the fixed time $T = n$. We shall make use of the i.i.d. inter-regeneration blocks as in the proof of Theorem 4 to show that

$$E_b(e^{-2\theta_b S_U + 2U\psi(\theta_b)} \mathbf{1}_{\{\theta_b S_n - n\psi(\theta_b) \geq n\phi(b)\}}) = O(n^{-1/2} e^{-2n\phi(b)}), \tag{A.9}$$

in which the additional $n^{-1/2}$ factor that is not present in Theorem 4 is due to the use of a local limit bound

$$\sup_{y \in \mathbf{R}, 1 \leq k \leq n^{1/2}} P_b\{y \leq \theta_b S_{n-k} - (n-k)\psi(\theta_b) < y + 1\} = O(n^{-1/2}), \tag{A.10}$$

in place of Blackwell's theorem for the renewal function (3.24). The proof of (A.10) is given in the next paragraph. Making use of Theorem 3 and Jensen's inequality and noting that $\{\theta_b S_n - n\psi(\theta_b) \geq n\phi(b)\} = \{S_n/n \geq b\}$, it can be shown by using arguments similar to those in the proof of Corollary 1 that

$$E_b \widehat{p}_n^2 / r^2(x_0; \theta_b) \leq E_b(e^{-2\theta_b S_U + 2U\psi(\theta_b)} \mathbf{1}_{\{\theta_b S_n - n\psi(\theta_b) \geq n\phi(b)\}}). \tag{A.11}$$

The desired conclusion then follows from (2.34) for the case $q = 0$ (see proof of Theorem 2 of [5]) together with (A.9) and (A.11).

Let $\widetilde{\xi}_i = \theta_b \xi_i - \psi(\theta_b)$ and $\widetilde{S}_i = \widetilde{\xi}_1 + \cdots + \widetilde{\xi}_i$. For the special case of i.i.d. nonlattice $\xi_i$ with variance $\sigma^2 > 0$, Theorem 1 of [19] yields

$$\begin{aligned} &P_b\{n\phi(b) - \sqrt{n}\sigma z - 1 \leq \widetilde{S}_n < n\phi(b) - \sqrt{n}\sigma z\} \\ &\qquad = (2\pi n\sigma^2)^{-1/2}[e^{-z^2/2} + o(1)] \qquad \text{as } n \to \infty, \end{aligned} \tag{A.12}$$

uniformly over $z \in \mathbf{R}$, so (A.10) holds. For Markov additive processes with nonlattice increments and satisfying the minorization condition (3.1), Chan and Lai [5] have shown that

$$\begin{aligned} &\int g(y) P_b\{X_n \in dy, n\phi(b) - \sqrt{n}\sigma z - 1 \leq \widetilde{S}_n < n\phi(b) - \sqrt{n}\sigma z\} \\ &\qquad = (2\pi n\sigma^2)^{-1/2} e^{-z^2/2} \left\{ \int g(y)\, d\pi(\theta_b) + o(1) \right\} \end{aligned}$$

for any nonnegative bounded measurable function $g$; see their (6.11) and the arguments leading to their (6.12). Taking $g \equiv 1$ then yields (A.12) and therefore also (A.10). When the Markov additive process has lattice increments, although one no longer has the precise formula (A.12), (A.10) still holds by a modification of these equations in [5], similar to that used to weaken (2.24) into (2.34) for lattice random walks.

To prove (A.9), let $\tau^*$ be the last regeneration time at or before time $n$. Then analogously to (3.25),

$$\begin{aligned}
&E_b(e^{2[n\phi(b)-\widetilde{S}_U]}\mathbf{1}_{\{\widetilde{S}_n\ge n\phi(b)\}})\\
&\quad=\sum_{k=0}^{n}\sum_{y=-\infty}^{\infty}E_b(e^{2[n\phi(b)-\widetilde{S}_U]}\mathbf{1}_{\{\widetilde{S}_n\ge n\phi(b)\}}\\
&\qquad\qquad\times\mathbf{1}_{\{\tau^*=n-k,n\phi(b)-y-1\le\widetilde{S}_{\tau^*}<n\phi(b)-y\}}).
\end{aligned}\tag{A.13}$$

Let $k\le n$. Using the decomposition $\widetilde{S}_U=\widetilde{S}_{n-k}+(\widetilde{S}_U-\widetilde{S}_{n-k})$ on the event $\{\tau^*=n-k\}$, we obtain

$$\begin{aligned}
&E_b(e^{2[n\phi(b)-\widetilde{S}_U]}\mathbf{1}_{\{\widetilde{S}_n\ge n\phi(b)\}}\mathbf{1}_{\{\tau^*=n-k,n\phi(b)-y-1\le\widetilde{S}_{\tau^*}<n\phi(b)-y\}})\\
&\quad=E_b(e^{2[n\phi(b)-\widetilde{S}_{n-k}]}\\
&\qquad\times\mathbf{1}_{\{\tau^*=n-k,n\phi(b)-y-1\le\widetilde{S}_{n-k}<n\phi(b)-y\}}e^{-2(\widetilde{S}_U-\widetilde{S}_{n-k})}\mathbf{1}_{\{\widetilde{S}_n\ge n\phi(b)\}})\\
&\quad\le e^{2(y+1)}E_b(e^{-2(\widetilde{S}_U-\widetilde{S}_{n-k})}\\
&\qquad\qquad\times\mathbf{1}_{\{\tau^*=n-k,\widetilde{S}_n-\widetilde{S}_{\tau^*}\ge y,n\phi(b)-y-1\le\widetilde{S}_{\tau^*}\le n\phi(b)-y\}}).
\end{aligned}\tag{A.14}$$

Conditioned on the event $\{\tau(i)=n-k,n\phi(b)-y-1\le\widetilde{S}_{\tau(i)}<n\phi(b)-y\}$, the vector $(\widetilde{S}_{\tau(i+1)}-\widetilde{S}_{\tau(i)},\widetilde{S}_n-\widetilde{S}_{\tau(i)},\tau(i+1)-\tau(i))$ has the same distribution as $(\widetilde{S}_\tau,\widetilde{S}_k,\tau)$ that is initialized at the regeneration distribution under $P_b$. Hence by (A.10), for $k\le n^{1/2}$,

$$\begin{aligned}
&E_b(e^{-2(\widetilde{S}_U-\widetilde{S}_{n-k})}\mathbf{1}_{\{\tau^*=n-k,\widetilde{S}_n-\widetilde{S}_{\tau^*}\ge y,n\phi(b)-y-1\le\widetilde{S}_{\tau^*}\le n\phi(b)-y\}})\\
&\quad=\sum_{i=0}^{\infty}E_b(e^{-2(\widetilde{S}_{\tau(i+1)}-\widetilde{S}_{\tau(i)})}\\
&\qquad\qquad\times\mathbf{1}_{\{\tau(i)=n-k,\tau(i+1)-\tau(i)>k,\widetilde{S}_n-\widetilde{S}_{\tau(i)}\ge y,n\phi(b)-y-1\le\widetilde{S}_{\tau(i)}<n\phi(b)-y\}})\\
&\quad=E_{\nu_b,b}(e^{-2\widetilde{S}_\tau}\mathbf{1}_{\{\widetilde{S}_k\ge y,\tau>k\}})\\
&\qquad\times\sum_{i=0}^{\infty}P_b\{\tau(i)=n-k,n\phi(b)-y-1\le\widetilde{S}_{\tau(i)}<n\phi(b)-y\}\qquad\text{(A.15)}\\
&\quad\le E_{\nu_b,b}(e^{-2\widetilde{S}_\tau}\mathbf{1}_{\{\widetilde{S}_k\ge y,\tau>k\}})\\
&\qquad\times P_b\{\tau(i)=n-k\text{ for some }i,\\
&\qquad\qquad n\phi(b)-y-1\le\widetilde{S}_{n-k}<n\phi(b)-y\}
\end{aligned}$$

$$= O(n^{-1/2}) E_{\nu_b,b}(e^{-2\widetilde{S}_\tau} \mathbf{1}_{\{\widetilde{S}_k \geq y, \tau > k\}}).$$

Moreover,

$$\begin{aligned}
&\sum_{0 \leq k \leq n^{1/2}} \sum_{y=-\infty}^{\infty} E_{\nu_b,b}(e^{2(y+1)-2\widetilde{S}_\tau} \mathbf{1}_{\{\widetilde{S}_k \geq y\}} \mathbf{1}_{\{\tau > k\}}) \\
&\quad \leq [e^2/(1-e^{-2})] \sum_{k=0}^{\infty} E_{\nu_b,b}[e^{2\widetilde{S}_k - 2\widetilde{S}_\tau} \mathbf{1}_{\{\tau > k\}}] \\
&\quad = [e^2/(1-e^{-2})] E_{\nu_b,b}\left[\sum_{k=0}^{\tau-1} e^{2\widetilde{S}_k - 2\widetilde{S}_\tau}\right] \\
&\quad = [e^2/(1-e^{-2})] E_\nu\left[\sum_{k=0}^{\tau-1} e^{2\widetilde{S}_k - \widetilde{S}_\tau}\right],
\end{aligned} \tag{A.16}$$

in which the last equality can be shown by using the same arguments as in (3.27). By (A.14)–(A.16),

$$\begin{aligned}
&\sum_{0 \leq k \leq n^{1/2}} \sum_{y=-\infty}^{\infty} E_b(e^{2[n\phi(b) - \widetilde{S}_U]} \mathbf{1}_{\{\widetilde{S}_n \geq n\phi(b)\}} \\
&\qquad \times \mathbf{1}_{\{\tau^* = n-k, n\phi(b) - y - 1 \leq \widetilde{S}_{\tau^*} < n\phi(b) - y\}}) \\
&\quad = O(n^{-1/2}) E_\nu\left(\sum_{k=0}^{\tau-1} e^{2\widetilde{S}_k - \widetilde{S}_\tau}\right) = O(n^{-1/2}),
\end{aligned} \tag{A.17}$$

where the last relation follows from Lemma 2 below.

Let $n^{1/2} < k < n$ and $\lambda = [4\psi(\theta_b) - \zeta]/2$. The bound in (A.15) can be modified to

$$\begin{aligned}
&E_b(e^{-2(\widetilde{S}_U - \widetilde{S}_{n-k})} \mathbf{1}_{\{\tau^* = n-k, \widetilde{S}_n - \widetilde{S}_{\tau^*} \geq y, n\phi(b) - y - 1 \leq \widetilde{S}_{\tau^*} < n\phi(b) - y\}}) \\
&\quad \leq E_{\nu_b,b}(e^{-2\widetilde{S}_\tau} \mathbf{1}_{\{\widetilde{S}_k \geq y, \tau > k\}}),
\end{aligned}$$

and therefore by (A.14), we can modify (A.16) and (A.17) to

$$\begin{aligned}
&\sum_{n^{1/2} < k < n} \sum_{y=-\infty}^{\infty} E_b(e^{2[n\phi(b) - \widetilde{S}_U]} \mathbf{1}_{\{\widetilde{S}_n \geq n\phi(b)\}} \\
&\qquad \times \mathbf{1}_{\{\tau^* = n-k, n\phi(b) - y - 1 \leq \widetilde{S}_{\tau^*} < n\phi(b) - y\}}) \\
&\quad \leq e^{-n^{1/2}\lambda} \sum_{n^{1/2} < k < n} e^{k\lambda} \sum_{y=-\infty}^{\infty} E_b(e^{2[n\phi(b) - \widetilde{S}_U]} \mathbf{1}_{\{\widetilde{S}_n \geq n\phi(b)\}}
\end{aligned} \tag{A.18}$$

$$\times \mathbf{1}_{\{\tau^*=n-k, n\phi(b)-y-1\le \widetilde{S}_{\tau^*}<n\phi(b)-y\}}\Big)$$

$$\le [e^{-n^{1/2}\lambda+2}/(1-e^{-2})]E_\nu\left(\sum_{k=0}^{\tau-1} e^{2\widetilde{S}_k-\widetilde{S}_\tau+k\lambda}\right) = O(e^{-\lambda\sqrt{n}})$$

since $E_\nu(\sum_{k=0}^{\tau-1} e^{2\widetilde{S}_k-\widetilde{S}_\tau+k\lambda})<\infty$ by Lemma 2 below. Finally, for the case $k=n$, by the Cauchy–Schwarz inequality,

$$\sum_{y=-\infty}^{\infty} E_b(e^{2[n\phi(b)-\widetilde{S}_U]}\mathbf{1}_{\{\widetilde{S}_n\ge n\phi(b)\}}\mathbf{1}_{\{\tau^*=0, n\phi(b)-y-1\le \widetilde{S}_0<n\phi(b)-y\}})$$

$$= E_b(e^{2[n\phi(b)-\widetilde{S}_\tau]}\mathbf{1}_{\{\widetilde{S}_n\ge n\phi(b),\tau>n\}}) \le e^{-n\lambda}E_b(e^{2[\widetilde{S}_n-\widetilde{S}_\tau]+n\lambda}\mathbf{1}_{\{\tau>n\}}) \tag{A.19}$$

$$= e^{-n\lambda}E(e^{2\widetilde{S}_n-\widetilde{S}_\tau+n\lambda}\mathbf{1}_{\{\tau>n\}})/r(x_0;\theta_b) = O(e^{-n\lambda}),$$

where the last relation follows from Lemma 2. From (A.13) and (A.17)–(A.19), (A.9) follows. □

LEMMA 2. *With the same notation and assumptions as in Theorem 6, let* $\widetilde{\xi}_i=\theta_b\xi_i-\psi(\theta_b)$, $\widetilde{S}_i=\widetilde{\xi}_1+\cdots+\widetilde{\xi}_i$ *and* $\lambda=[4\psi(\theta_b)-\zeta]/2$. *Then*

$$E_\nu\left(\sum_{k=0}^{\tau-1} e^{2\widetilde{S}_k-\widetilde{S}_\tau}\right)+E_\nu\left(\sum_{k=0}^{\tau-1} e^{2\widetilde{S}_k-\widetilde{S}_\tau+k\lambda}\right)<\infty,$$

$$E(e^{2\tilde{S}_n-\tilde{S}_\tau+n\lambda}\mathbf{1}_{\{\tau>n\}})=O(1),$$

*where* $\tau$ *is the regeneration time as in* (A.15)–(A.19).

PROOF. Let $W_k=e^{2\widetilde{S}_k+k\lambda}\mathbf{1}_{\{\tau>k\}}$ and $Y=e^{-\widetilde{S}_\tau}$. Then

$$E_\nu\left(\sum_{k=0}^{\tau-1} e^{2\widetilde{S}_k-\widetilde{S}_\tau+k\lambda}\right) = E_\nu\left(\sum_{k=0}^{\infty} W_kY\right) = E_\nu\left(\sum_{k=0}^{\infty} W_kY\mathbf{1}_{\{\tau>k\}}\right)$$

$$\le \tfrac{1}{2}\left[E_\nu\left(\sum_{k=0}^{\infty} W_k^2\right)+E_\nu(\tau Y^2)\right]$$

$$=\{\ell_\nu(4\theta_b,\zeta)+E_\nu(\tau Y^2)\}/2.$$

Since $(-2\theta_b, 2\psi(\theta_b))\in\Omega$ and $\Omega$ is open, $E_\nu(\tau Y^2)<\infty$. Therefore

$$E_\nu\left(\sum_{k=0}^{\tau-1} e^{2\tilde{S}_k-\tilde{S}_\tau+k\lambda}\right)<\infty.$$

Since $\lambda>0$, this implies that $E_\nu(\sum_{k=0}^{\tau-1} e^{2\tilde{S}_k-\tilde{S}_\tau})$ is also finite. By the Cauchy–Schwarz inequality,

$$E(e^{2\tilde{S}_n-\tilde{S}_\tau+n\lambda}\mathbf{1}_{\{\tau>n\}}) \le [Ee^{4\tilde{S}_n+2n\lambda}\mathbf{1}_{\{\tau>n\}}]^{1/2}[Ee^{-2\tilde{S}_\tau}]^{1/2}=O(1),$$

since $Ee^{-2\tilde{S}_\tau} = EY^2 < \infty$ and since $E(e^{4\widetilde{S}_n + 2n\lambda}\mathbf{1}_{\{\tau>n\}}) = E(e^{4\theta_b S_n - n\zeta}\mathbf{1}_{\{\tau>n\}}) \le \ell_{x_0}(4\theta_b, \zeta) < \infty$. □

Proof of Theorem 7. The first step is to generalize Theorem 2 of [4] to the Markov case. This can be done by combining the basic ideas of the proof of that theorem with those of the proof of Theorem 3 of [5]. Assuming $r = 1$ without loss of generality, we can use these arguments to show that (2.24) still holds in the nonlattice case and that (2.34) holds without the nonlattice assumption. Whereas $p_n = E_{Q_n^*}\widehat{p}_n$ can be analyzed by using Chan and Lai's [5] truncation argument to handle $1/r(X_n^{(i)}; \theta_\mu)$ in (2.15) [see (3.34)], the analysis of $E_{Q_n^*}\widehat{p}_n^2$ involves $1/r^2(X_n^{(i)}; \theta_\mu)$ which does not relate to the finiteness of eigenmeasures via the truncation argument. For the special case $d = 1$ and $g(\mu) = \mu$, the proof of Theorem 6 uses regeneration and Theorem 3 to circumvent this difficulty. Note that in this special case, the exponential tilting involves a single $\theta_b$ instead of a mixture of $\theta_\mu$'s. We can use geometric integration over a suitably chosen tubular neighborhood of the manifold $M$ as in the proof of Theorem 2 of [4] to piece together the conclusions of Theorem 6 for the local tiltings. □

Proof of Corollary 4. Let $Q_n = \sum_{\mu \in G} \omega_{\mu,n} P_{\mu,n}$. Arguments similar to those in (4.6) can be used to show

$$E_{Q_n}[L_n^2 \mathbf{1}_{\{g(S_n/n) \ge b\}}] \le \sum_{\mu \in G} \omega_{\mu,n}^{-1} E_\mu\left[\left(\frac{dP_n}{dP_{\mu,n}}\right)^2 \mathbf{1}_{\{S_n/n \in H(\mu)\}}\right]. \tag{A.20}$$

Noting that $\{S_n/n \in H(\mu)\} = \{\theta_\mu' S_n - n\psi(\theta_\mu) \ge n\phi(\mu)\}$, it follows from the proof of Theorem 6 [in particular from the multidimensional versions of (A.9) and (A.11)] that

$$E_\mu\left[\left(\frac{dP_n}{dP_{\mu,n}}\right)^2 \mathbf{1}_{\{S_n/n \in H(\mu)\}}\right] = O(n^{-1/2} e^{-2n\phi(\mu)}), \qquad \mu \in G. \tag{A.21}$$

By (2.34) with $q = 0$ and with $c$ replaced by the more general $c/r$, $p_n$ is of the order $n^{-1/2} e^{-bn/r}$, and therefore $Q_n$ is asymptotically optimal in view of (4.8), (A.20) and (A.21). □

## REFERENCES

[1] Athreya, K. B. and Ney, P. (1978). A new approach to the limit theory of recurrent Markov chains. *Trans. Amer. Math. Soc.* **245** 493–501. MR0511425

[2] Bucklew, J. A., Ney, P. and Sadowsky, J. S. (1990). Monte Carlo simulation and large deviation theory for uniformly recurrent Markov chains. *J. Appl. Probab.* **27** 44–59. MR1039183

[3] Bucklew, J. A., Nitinawarat, S. and Wierer, J. (2004). Universal simulation distributions. *IEEE Trans. Inform. Theory* **50** 2674–2685. MR2096986

[4] Chan, H. P. and Lai, T. L. (2000). Asymptotic approximations for error probabilities of sequential or fixed sample size tests in exponential families. *Ann. Statist.* **28** 1638–1669. MR1835035

[5] Chan, H. P. and Lai, T. L. (2003). Saddlepoint approximations and nonlinear boundary crossing probabilities of Markov random walks. *Ann. Appl. Probab.* **13** 395–429. MR1970269

[6] Chan, H. P. and Lai, T. L. (2005). Importance sampling for generalized likelihood ratio procedures in sequential analysis. *Sequential Anal.* **24** 259–278. MR2187338

[7] Collamore, J. F. (2002). Importance sampling techniques for the multidimensional ruin problem for general Markov additive sequences of random vectors. *Ann. Appl. Probab.* **12** 382–421. MR1890070

[8] Dupuis, P. and Wang, H. (2005). Dynamic importance sampling for uniformly recurrent Markov chains. *Ann. Appl. Probab.* **15** 1–38. MR2115034

[9] Feller, W. (1971). *An Introduction to Probability Theory and Its Applications* **2**, 2nd ed. Wiley, New York.

[10] Glasserman, P. and Wang, Y. (1997). Counterexamples in importance sampling for large deviation probabilities. *Ann. Appl. Probab.* **7** 731–746. MR1459268

[11] Iscoe, I., Ney, P. and Nummelin, E. (1985). Large deviations of uniformly recurrent Markov additive processes. *Adv. in Appl. Math.* **6** 373–412. MR0826590

[12] Lehtonen, T. and Nyrhinen, H. (1992). Simulating level-crossing probabilities by importance sampling. *Adv. in Appl. Probab.* **24** 858–874. MR1188956

[13] Lehtonen, T. and Nyrhinen, H. (1992). On asymptotically efficient simulation of ruin probabilities in a Markovian environment. *Scand. Actuar. J.* **1** 60–75. MR1193671

[14] Meyn, S. P. and Tweedie, R. L. (1993). *Markov Chains and Stochastic Stability.* Springer, London. MR1287609

[15] Ney, P. and Nummelin, E. (1987). Markov additive processes. I. Eigenvalues properties and limit theorems. *Ann. Probab.* **15** 561–592. MR0885131

[16] Ney, P. and Nummelin, E. (1987). Markov additive processes. II. Large deviations. *Ann. Probab.* **15** 593–609. MR0885132

[17] Sadowsky, J. S. and Bucklew, J. A. (1990). On large deviations theory and asymptotically efficient Monte Carlo estimation. *IEEE Trans. Inform. Theory* **36** 579–588. MR1053850

[18] Siegmund, D. (1976). Importance sampling in the Monte Carlo study of sequential tests. *Ann. Statist.* **4** 673–684. MR0418369

[19] Stone, C.(1965). Local limit theorem for nonlattice multi-dimensional distribution functions. *Ann. Math. Statist.* **36** 546–551. MR0175166

Department of Statistics
and Applied Probability
National University of Singapore
Singapore 119260
Republic of Singapore
E-mail: stachp@nus.edu.sg

Department of Statistics
Stanford University
Stanford, California 94305-4065
USA
E-mail: lait@stat.stanford.edu